\documentclass[11pt]{amsart}
\usepackage[margin=2cm]{geometry}
\usepackage[english]{babel}
\usepackage{mathrsfs}
\usepackage[pagebackref]{hyperref}
\usepackage{enumitem}
\usepackage{graphicx,color}
\usepackage{comment}
\usepackage{url}
\usepackage{epstopdf}
\usepackage{float}
\usepackage{caption}
\usepackage{algorithm}
\usepackage{algpseudocode}
\usepackage{multicol}
\usepackage{tikz}
\usetikzlibrary{calc} 
\usepackage{xcolor}
 \definecolor{light-gray}{gray}{0.9}

\setenumerate[1]{label=(\roman*), ref=(\roman*)}
\usepackage[all]{xy}% Para los diagramas.
\newtheorem{thm}{Theorem}[section]
\newtheorem{prop}[thm]{Proposition}
\newtheorem{definition}[thm]{Definition}

\newcommand{\Flder}{\rightarrow}
\newcommand{\der}{\partial}
\newtheoremstyle{obs}% name
  {3pt}%      Space above
  {3pt}%      Space below
  {}%         Body font
  {}%         Indent amount (empty = no indent, \parindent = para indent)
  {\bfseries}% Thm head font
  {.}%        Punctuation after thm head
  {.5em}%     Space after thm head: " " = normal interword space;
  {}%         Thm head spec (can be left empty, meaning `normal')
\theoremstyle{obs}
\newtheorem{remark}[thm]{Remark}

\newtheorem*{ex*}{Example}

\newcommand{\R}{\mathbb{R}}      %Numeros reales
\newcommand{\lp}{\left(}
\newcommand{\rp}{\right)}
\newcommand{\lc}{\left\{}
\newcommand{\rc}{\right\}}

\def\qed{\ifvmode\removelastskip\fi
{\unskip\nobreak\hfil\penalty50\hbox{}\nobreak\hfil \hbox{\vrule
height1.2ex width1.2ex}\parfillskip=0pt \finalhyphendemerits=0
\par \smallskip}}
\begin{document}

\title[Nonholonomic Energy-preserving integrators]{\sc Energy-Preserving Integrators Applied to Nonholonomic Systems}

\author[E. Celledoni,  M. Farr\'e, E. H{\o}iseth, D. Mart{\'\i}n de Diego]{Elena Celledoni$^1$, Marta Farr\'e Puiggal{\'\i}$^2$, Eirik Hoel H{\o}iseth$^1$\\  David Mart{\'\i}n de Diego$^2$}

\address{$^1$ 	Department of Mathematical Sciences, NTNU, 
7491 Trondheim,	Norway \href{elenac@math.ntnu.no}{} \href{eirik.hoiseth@math.ntnu.no}{}
}
\email[Elena Celledoni]{elenac@math.ntnu.no} \email[Eirik Hoel H{\o}iseth]{eirik.hoiseth@math.ntnu.no}
\address{
$^2$  Instituto de Ciencias Matem\'aticas (CSIC-UAM-UC3M-UCM),
   Calle Nicol\'as Cabrera 15, Campus UAM, Cantoblanco, 
   Madrid, 28049, Spain 	\href{mfarrepu7@gmail.com}{} \href{david.martin@icmat.es}{}}
\email[Marta Farr\'e Puiggal{\'\i}]{marta.farre@icmat.es} \email[David Mart{\'\i}n de Diego]{david.martin@icmat.es}

\begin{abstract}
We introduce energy-preserving integrators for
nonholonomic mechanical systems.  We will see that the nonholonomic  dynamics is completely determined  by a triple $({\mathcal D}^*, \Pi, \mathcal{H})$, where ${\mathcal D}^*$ is the dual of the vector bundle determined by the nonholonomic constraints, $\Pi$ is an almost-Poisson bracket (the nonholonomic bracket) and $ \mathcal{H}: {\mathcal D}^*\rightarrow \R$ is a Hamiltonian function. 
For this triple, we can  apply energy-preserving integrators, in particular,  we show that discrete gradients can be used in the numerical integration of nonholonomic dynamics.
By construction, we achieve preservation of the constraints and of the energy of the nonholonomic system. Moreover, to facilitate their applicability to complex systems which cannot be easily transformed into the aforementioned almost-Poisson form, we rewrite our integrators using just the initial information of the nonholonomic system. 
The derived procedures are tested on several examples: A chaotic quartic nonholonomic mechanical system, the Chaplygin sleigh system, the Suslov problem and a continuous gearbox driven by an asymmetric pendulum. Their performace is compared with other standard methods in nonholonomic dynamics, and their merits verified in practice.

\end{abstract}

\maketitle

\section{Introduction}

Geometric integrators are numerical methods for differential equations which preserve structural properties such as constants of the motion, symplectic or Poisson structures, phase-space volume,  different symmetries of the system or isospectrality. Preservation of structural properties is often desirable to achieve correct qualitative behaviour and long time stability \cite{hairer,mcqu,serna}.

In this paper, we address the construction of geometric integrators for mechanical systems subjected to nonholonomic constraints. 
There is considerable interest in the study of nonholonomic systems since nonholonomic constraints are present in a great variety of mechanical systems in engineering and robotics. For instance, they describe the dynamics of wheeled vehicles, manipulation devices and locomotion systems (see \cite{bloch,BKMM96,cortes02,cortmart,Murray,neimarkfufaev} and references therein).  

In the unconstrained case, or when the constraints are holonomic, mechanical systems have many distinguishing geometric features. Among the most important are the preservation of energy, the symplectic form constructed from the Lagrangian (Poincar\'e-Cartan 2-form) and the momentum map in the presence of symmetries according to the Noether theorem. As we will see, when we are dealing with nonholonomic constraints this symplectic form is no longer preserved, and the momentum map is not in general conserved in the presence of symmetries. However, the energy is still a conservation law for the system in the case of linear constraints. We therefore focus our attention on the exact preservation of energy, using geometric integrators, while writing the equations of motion in a format which ensures the nonholonomic constraints are satisfied. 

The proposed approach is different from other recent approaches  such as \cite{cortes02,dLdDSM04,feze,FeIgMa,KoMaFe,perlmutter06} where the authors have introduced numerical integrators for nonholonomic systems with very good energy behavior, and properties such as the
preservation of the discrete nonholonomic momentum map.

The ingredients of the theory of nonholonomic mechanics are a Lagrangian function $L: TQ\rightarrow \R$ of mechanical type, that is, kinetic minus potential energy, and a vector subbundle $\tau_{\mathcal D}: {\mathcal D}\rightarrow Q$ of $TQ$. This vector subbundle determines the nonholonomic constraints, as we will see in the next section. To develop integrators we first introduce a Hamiltonian description of nonholonomic mechanics in terms of an almost-Poisson bracket. Using the Riemannian metric determined by the kinetic energy, and the standard symplectic structure on $T^*Q$, we can induce a linear almost-Poisson structure $\Pi$ on the dual bundle $\pi_{D}: {\mathcal D}^*\rightarrow Q$. This so-called nonholonomic bracket is isomorphic to the nonholonomic bracket considered in \cite{SM1994}.

Now, working on the ``Hamiltonian system" determined by the triple given by (i) ${\mathcal D}^*$ as new phase space, (ii) the  almost-Poisson bracket $\Pi$ and (iii) the induced Hamiltonian function ${\mathcal H}:  {\mathcal D}^*\rightarrow {\mathbb R}$, we apply energy-preserving integrators to simulate its dynamics. 
This is a coherent approach since the unique generic quantity preserved by the flow of the corresponding Hamiltonian vector field to ${\mathcal H}$ is precisely the Hamiltonian function. The resulting integrators preserve by construction both the energy and nonholonomic constraints. 

To approximate the solution while preserving the energy of the initial nonholonomic problem we utilize a class of geometric integrators called discrete gradient methods. Consider an ODE which can be written in skew-gradient form, i.e. $\dot{x}=\Pi(x)\nabla {\mathcal H}(x)$ with $x\in \mathbb{R}^N$ and $\Pi(x)$ a skew-symmetric matrix. 
In \cite{McQuRo} it is shown that any ODE with a  generic first integral ${\mathcal H}$ can be put into skew-gradient form. For a generalisation of these ideas to the case where the configuration space is a Lie group or a homogeneous manifold see \cite{CeOw}. 

Now, discrete gradient methods are based on the following construction: Let $x \approx x(nh)$ and $x' \approx x((n+1)h)$. Using a discrete gradient $\bar{\nabla}{\mathcal H}(x, x')$, which is an appropriate approximation of the gradient of ${\mathcal H}$ (see Section \ref{section3} for details), it is possible to define a class of integrators 
\[
\frac{x'-x}{h}=\tilde{\Pi}(x, x')\bar{\nabla}H(x, x') \, ,
\]
which preserve the first integral ${\mathcal H}$ exactly, i.e. ${\mathcal H}(x)={\mathcal H}(x')$.
Here $\tilde{\Pi}(x, x')$ is a skew-symmetric matrix approximating $\Pi(x)$. It can be shown that, in $\R^n$, any first integral-preserving (direct) integrator can be written as a discrete gradient method \cite{GONZ,Norton,QC1996,QT1996}. 

For a given nonholonomic mechanical system, the equations of motion in canonical coordinates are generally assumed known. A potential obstacle in applying a discrete gradient method directly to the adapted coordinate system is the need for the user to analytically derive these equations. For this reason, we propose a reformulation of the methods using just information from the original system in canonical coordinates. With this approach the analytic reformulation of the system in adapted coordinates is avoided. 
 
  The outline of the paper is as follows. In the next section we will describe the geometric framework for nonholonomic mechanics. The main objective is describing its dynamics as a Hamiltonian system on a vector bundle equipped with an almost-Poisson bracket. The resulting equations of motion in adapted coordinates are seen to be explicitly given in skew-gradient form.
 In Section \ref{section3} we apply discrete gradient integrators to the derived formulation to get energy-preserving integrators for nonholonomic systems. We then rewrite these integrators in an equivalent form by using only the information from the original nonholonomic system.
 Finally, in Section \ref{section4}, we verify the properties and the performance of our integration techniques, applying them to several interesting examples: 
The chaotic quartic nonholonomic mechanical system, the Chaplygin sleigh system, the Suslov problem and the continuous gearbox driven by an asymmetric pendulum. Our methods are compared with other well known numerical methods for nonholonomic mechanics.

\section{Nonholonomic systems}\label{section2}
A nonholonomic system is a mechanical system with external constraints on the velocities \cite{bloch,BKMM96,cortes02}. 
We will only consider linear velocity constraints, since this is the case in most examples. Linear velocity constraints are constraints that are specified by a regular $C^{\infty}$-distribution ${\mathcal D}$  on the configuration manifold $Q$, or equivalently, by a vector subbundle $\tau_{\mathcal D}: {\mathcal D}\rightarrow Q$ of $TQ$ with canonical inclusion $i_{\mathcal D}: {\mathcal D}\hookrightarrow TQ$. Therefore, we will say that a curve $\gamma: I\subseteq {\mathbb R}\rightarrow Q$ satisfies the constraints given by ${\mathcal D}$ if 
\begin{equation}\label{curve}
\gamma'(t)=\frac{d\gamma}{dt}(t)\in {\mathcal D}_{\gamma(t)} \hbox{    for all   } t\in I \, . 
\end{equation}
We say that ${\mathcal D}$ is holonomic if ${\mathcal D}$ is integrable or involutive, that is, for any vector fields  $X, Y\in {\mathfrak X}(Q)$ taking values on ${\mathcal D}$, it holds that the vector field $[X, Y]$ also takes values pn ${\mathcal D}$. A regular linear velocity constraint ${\mathcal D}$ is nonholonomic if  it is not holonomic. 
Observe that in the case of holonomic constraints all the curves through a point $q\in Q$ satisfying the constraints must lie on the maximal integral manifold for ${\mathcal D}$ through $q$. 

Let $\dim Q=n$. Locally if $(q^i)$, $1\leq i\leq n$ are coordinates on $Q$ and 
$(q^i, \dot{q}^i)$ are the induced coordinates on $TQ$,  the linear constraints are written as 
\begin{equation*}\label{LC}
\mu^{\alpha}_{i}\lp q\rp\dot q^{i}=0,\hspace{2mm} m+1\leq \alpha\leq n\, ,
\end{equation*}
where $\mbox{rank}\lp\mathcal{D}\rp=m\leq n$. The annihilator $\mathcal{D}^{\circ}$ is locally given by
\[
\mathcal{D}^{\circ}=\mbox{span}\lc\mu^{\alpha}=\mu_{i}^{\alpha}(q)\,dq^{i};\hspace{1mm} m+1\leq \alpha\leq n\rc\, ,
\]
where the 1-forms $\mu^{\alpha}$ are independent.
Equivalently, we can find independent vector fields $\{X_a\}$, $1\leq a\leq m$ such that 
\[
{\mathcal D}_q=\hbox{span}\{X_a\}\, .
\]
Observe that $\mu^{\alpha}(X_a)=0$, for all $m+1\leq \alpha\leq n$ and $1\leq a\leq m$. 

\medskip

\hspace{0.015\linewidth}
\begin{minipage}{0.95\textwidth}
\begin{ex*}[Rolling disk]
{	\rm \small
One of the simplest examples of a nonholonomic system is the unicycle, that is, a disk of radius $r$ which rolls on a horizontal plane, and always 
remains exactly upright, see for instance \cite{bloch}.  The coordinates $(x_1, x_2, \theta, \phi)$ describe the possible configurations of the system, where $(x_1, x_2)$ are the coordinates of the contact point with the $x_1x_2$-plane,  $\theta$ the heading angle and $\phi$  the self-rotation angle. The system is shown in Figure \ref{fig:disk}.
Not all the velocities are admissible for this system , since the constraint that the disk roll without slipping is specified by the  linear velocity constraints
\begin{equation}
\label{diskConst}
\dot{x}_1-r \dot{\phi}\cos\theta=0\, ,\qquad \dot{x}_2-r \dot{\phi}\sin\theta=0\, ,
\end{equation}
and therefore, 
\[
{\mathcal D}=\hbox {span} \left\{X_1=r\cos\theta \frac{\partial}{\partial x_1}+r\sin \theta \frac{\partial}{\partial x_2}+\frac{\partial}{\partial \phi}\, , X_2=\frac{\partial}{\partial \theta}  \right\} \, .
\]
The constraints are nonholonomic since the distribution ${\mathcal D}$ is not involutive. Indeed,
\[
[X_1, X_2]=r\sin\theta\frac{\partial}{\partial x_1}-r\cos\theta\frac{\partial}{\partial x_2}\, ,
\]
and $[X_1, X_2](q)\notin {\mathcal D}_q$ for all $q\in Q$.
}

\begin{center}
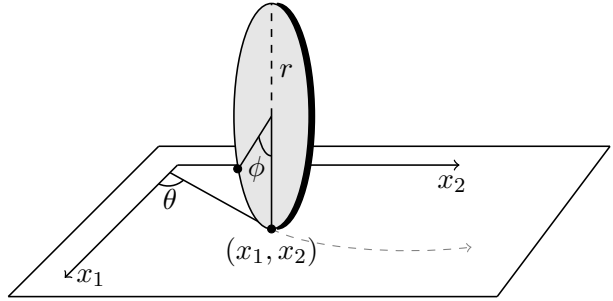


\begin{tikzpicture}
    
   % Surface
	\draw[semithick] (0,0) -- (2,2);
	\draw[semithick] (2,2) -- (8,2);
	\draw[semithick] (8,2) -- (6.5,0);
	\draw[semithick] (6.5,0) -- (0,0);
	
	% Line to disk
	\draw[semithick] (2.15,1.65) -- (3.5,0.9);
	
	% Axes
	\draw[->,semithick] (2.25,1.75) -- (6,1.75);
	\draw[->,semithick] (2.25,1.75) -- (0.75,0.25);	
	%\draw[->,semithick] (2.25,1.75) -- (2.25,4.75);
	
	% theta arc
	\draw[semithick] (2,1.5) arc (-115:-47:0.3 and 0.3);
				
	% rolling disk		
	\draw[semithick,fill=black] (3.57,0.9) arc (-90:90:0.5 and 1.5);
	\draw[semithick,fill=light-gray] (3.5,2.4) ellipse (0.5 and 1.5);
	
   % Trajectory	
	\draw[->,dashed,color=gray] (3.5,0.9) arc (-155:-65:2 and 0.5);
	
	% phi lines and arcs
	\draw[semithick] (3.5,0.9) -- (3.5,2.4);
	\draw[semithick] (3.5,2.4) -- (3.06,1.7);
	\draw[semithick] (3.33,2.13) arc (-145:-90:0.2 and 0.6);
	
	% Radius line
	\draw[dashed,semithick] (3.5,2.4) -- (3.5,3.9);

    % Points on unicycle
    \draw [fill=black] (3.5,0.9) circle[radius=1.5pt];
    %\draw [fill=black] (3.5,2.4) circle[radius=1.5pt];
    \draw [fill=black] (3.05,1.7) circle[radius=1.5pt];
    
    % Nodes
    	\draw (3.5,0.6) node {$(x_1,x_2)$};
    \draw (3.7,3) node {$r$};
	\draw (5.9,1.5) node {$x_2$};
	\draw (1.1,0.25) node {$x_1$};
	%\draw (2,4.75) node {$z$};
	\draw (2.15,1.3) node {$\theta$};
	\draw (3.3,1.7) node {$\phi$};	
\end{tikzpicture}
\captionof{figure}{The geometry of the rolling disk.}
\label{fig:disk}
\end{center}

\hfill $\Diamond$
\end{ex*}
\end{minipage}

\subsection{Lagrangian equations for nonholonomic systems}
In addition to the constraints, the dynamics is specified by a Lagrangian function $L: TQ\rightarrow {\mathbb R}$, assumed to be of mechanical type, that is,
\begin{equation*}\label{Mec}
L\lp v_{q}\rp=\frac{1}{2}g\lp v_{q},v_{q}\rp-V\lp q\rp\, ,\quad v_{q}\in T_{q}Q \, ,
\end{equation*}
where $g$ is a Riemannian metric on the configuration space $Q$ and $V: Q\rightarrow {\mathbb R}$ a potential function.  The Lagrangian is written in coordinates $(q^i, \dot{q}^i)$ as 
\[
L(q^i, \dot{q}^i)=\frac{1}{2}\textbf{\textit{g}}_{ij}(q)\dot{q}^i\dot{q}^j-V(q)\, ,
\]
where $\textbf{\textit{g}}_{ij}=g\lp\der/\der q^{i},\der/\der q^{j}\rp$, $1\leq i,j\leq n$.

\medskip

\hspace{0.015\linewidth}
\begin{minipage}{0.95\textwidth}
\begin{ex*}[Rolling disk, continued]
{	\rm \small
For the case of the rolling disk the Riemannian metric is 
\[
g=m\,dx_1\otimes dx_1 +m\,dx_2\otimes dx_2 +J_{\theta}\, d\theta\otimes d\theta+J_{\phi}\,d\phi\otimes d\phi \, ,
\]
where $m$ is the mass of the disk and $J_{\theta}$ and $J_{\phi}$ are the moment of inertia about the $\theta$ and $\phi$ axis respectively. We assume that the Lagrangian is purely kinetic, i.e. $V=0$, and thus
\[
L(x_1, x_2, \theta, \phi, \dot{x}_1, \dot{x}_2, \dot{\theta}, \dot{\phi})=\frac{1}{2}\left(m\,\dot{x}_1^2 +m\,\dot{x}_2^2+J_{\theta}\, \dot{\theta}^2+J_{\phi}\,\dot{\phi}^2\right)\, .   
\]\hfill $\Diamond$
}
\end{ex*}
\end{minipage}

\medskip

In nonholonomic mechanics, the  equations of motion are completely determined by the   La\-gran\-ge-d'Alembert principle. This principle states that a curve $q:I\subset \R\Flder Q$  is an admissible motion of the system if
\[
\delta\mathcal{J}=\delta\int^{T}_{0}L\lp q\lp t\rp, \dot q\lp t\rp\rp dt=0\, ,
\]
for all variations such that $\delta q\lp t\rp\in\mathcal{D}_{q\lp t\rp}$, $0\leq t\leq T$, $\delta q\lp 0\rp=\delta q\lp T\rp=0$. The velocity of the curve itself must also satisfy the constraints according to \eqref{curve}.  From the Lagrange-d'Alembert principle, we arrive at the well-known {\bf nonholonomic equations}
\begin{subequations}\label{LdAeqs}
\begin{align}
\frac{d}{dt}\lp\frac{\der L}{\der\dot q^{i}}\rp-\frac{\der L}{\der q^{i}}&=\lambda_{\alpha}\mu^{\alpha}_{i}\, ,\label{Con-1}\\\label{Con-2}
\mu_{i}^{\alpha}(q)\,\dot q^{i}&=0\, ,
\end{align}
\end{subequations}
where $\lambda_{\alpha}$, $m+1\leq \alpha \leq n$, is a set of Lagrange multipliers. The right-hand side of Equation (\ref{Con-1}) represents the force induced by the constraints, while Equation (\ref{Con-2}) gives the linear velocity constraints themselves.

It is important to stress that in Equations (\ref{LdAeqs}) it is necessary to use the Lagrangian defined on the full space $TQ$ instead of working with the restriction of $L$ to ${\mathcal D}$ (where we now consider ${\mathcal D}$ as a vector subbundle of $TQ$). Applying standard variational techniques and using $l=L|_{{\mathcal D}}$ we would derive a different set of equations than (\ref{LdAeqs}), which are not valid for nonholonomic mechanics. These other equations are called vakonomic equations, or variational constrained equations in the literature, see for instance \cite{arnoldIII}.   

\subsubsection{Adapted coordinates}
Equations (\ref{LdAeqs}) are derived using a set of coordinates $(q^i)$ on $Q$, and the induced coordinates $(q^i, \dot{q}^i)$ on $TQ$ by the canonical coordinate frame 
$
\left\{
\frac{\partial}{\partial q^i}
\right\}
$, $1\leq i\leq n$.
Any element $v_q\in T_qQ$ can therefore be written univocally as
\[
v_q=\dot{q}^i\frac{\partial}{\partial q^i}\Big|_q\, .
\]
In the case of nonholonomic mechanics it can be useful to adapt the chosen frame to the linear velocity constraints. Specifically we consider a basis of vector fields 
$\{X_a, X_{\alpha}\}$, $1\leq a\leq m$ and $m+1\leq \alpha\leq n$, such that locally 
\[
{\mathcal D}_q=\hbox{span}\{X_a(q)\} \quad \hbox{and}\quad {\mathcal D}^{\perp, g}_q=\hbox{span}\{X_{\alpha}(q)\}\, ,
\]
where ${\mathcal D}^{\perp, g}_q$ is the (Riemannian)-orthogonal to ${\mathcal D}$, i.e.
\[
g(X_a, X_{\alpha})=0\, , \quad 1\leq a\leq m\, \quad \hbox{and} \quad m+1\leq \alpha\leq n\, .
\]
Observe that $T_qQ={\mathcal D}_q\oplus {\mathcal D}^{\perp, g}_q$.

The adapted  basis $\{X_a, X_{\alpha}\}$  induces  a new  set of coordinates on the tangent bundle $(q^i, y^a, y^{\alpha})$ (also called quasi-velocities)
so that now
\[
v_q=y^a X_a(q)+y^{\alpha} X_{\alpha}(q)\, .
\]
Observe that the elements $v_q\in {\mathcal D}_q$ are distinguished by $y^{\alpha}=0$. Therefore $y^{\alpha}=0$ expresses the nonholonomic constraints in the adapted basis. Consequently ${\mathcal D}$ is completely described by coordinates $(q^i, y^a)$.

\medskip

\hspace{0.015\linewidth}
\begin{minipage}{0.95\textwidth}
\begin{ex*}[Rolling disk, continued]
{	\rm \small  
We take an adapted basis $\{X_1, X_2, X_3, X_4\}$, where
\begin{align*}
{\mathcal D}&=\hbox {span} \left\{X_1=r\cos\theta \frac{\partial}{\partial x_1}+r\sin \theta \frac{\partial}{\partial x_2}+\frac{\partial}{\partial \phi}, X_2=\frac{\partial}{\partial \theta}  \right\} \, ,\\
{\mathcal D}^{\perp, g}&=\hbox {span} \left\{
X_3=\frac{1}{m}\frac{\partial}{\partial x_1}-\frac{r}{J_{\phi}}\cos\phi \frac{\partial}{\partial \phi}, X_4=\frac{1}{m}\frac{\partial}{\partial x_2}-\frac{r}{J_{\phi}}\sin\phi \frac{\partial}{\partial \phi}  \right\} \, .\\
\end{align*}
This induces coordinates $(x_1, x_2, \theta, \phi, y^1, y^2, y^3, y^4)$ on $TQ$, which are related to the standard coordinates as follows: 
\begin{align*}
\dot{x}_1&=ry^1\cos\theta+\frac{y^3}{m}  \, ,   \\
\dot{x}_2&= ry^1\sin\theta +\frac{y^4}{m} \, ,  \\
\dot{\theta}&=y^2  \, , \\
\dot{\phi}&=y^1-\frac{r}{J_{\phi}}y^3\cos\phi-\frac{r}{J_{\phi}}y^4\sin\phi \, .
\end{align*}
Observe that the linear constraints have the simple form  $y^3=0, y^4=0$ in the adapted basis. \hfill $\Diamond$
}
\end{ex*}
\end{minipage}

\medskip

\subsubsection{Equations of motion in adapted coordinates}

We now want to rewrite the equations of motion of the nonholonomic system in terms of the coordinates $(q^i, y^a, y^{\alpha})$, instead of the canonical coordinates $(q^i, \dot{q}^i)$. 
Consider first Equation (\ref{Con-1}). We can split it as the following system of equations: 
\begin{subequations}
\begin{align}
0&=X^i_a\frac{d}{dt}\lp\frac{\der L}{\der\dot q^{i}}\rp-X^i_a\frac{\der L}{\der q^{i}} -\lambda_{\beta}\mu^{\beta}_{i}X^i_a=X^i_a\frac{d}{dt}\lp\frac{\der L}{\der\dot q^{i}}\rp-X^i_a\frac{\der L}{\der q^{i}}  \, ,  \label{equation1} \\
0&=X^i_{\alpha}\frac{d}{dt}\lp\frac{\der L}{\der\dot q^{i}}\rp-X^i_{\alpha}\frac{\der L}{\der q^{i}} -\lambda_{\beta}\mu^{\beta}_{i}X^i_{\alpha} \, ,  \label{equation2}
\end{align}
\end{subequations}
with $1 \leq a \leq m$, $m+1\leq \alpha, \beta \leq n$. In Equation (\ref{equation1}) we have used that $X_a(q)\in {\mathcal D}_q$. Observe that Equation (\ref{equation2}) uniquely gives information about the value of the Lagrange multipliers since $g$ is a Riemannian metric and therefore $(\mu^{\beta}_{i}X^i_{\alpha})$ is a regular matrix. As we are not interested in the evolution of the Lagrange multipliers we discard this second set of equations. 

Define the Lagrangian in adapted coordinates as $\tilde{L}(q^i, y^a, y^{\alpha}) :=L(q^i, X^i_a y^a+X^i_{\alpha}y^{\alpha})$. We want to express Equation (\ref{equation1}) in terms of $\tilde{L}$. To this end observe that 
\begin{align*}
\frac{\partial \tilde{L}}{\partial y^a}&=X^i_a\frac{\der L}{\der \dot{q}^{i}} \, , \\
\frac{\partial \tilde{L}}{\partial y^{\alpha}}&=X^i_{\alpha}\frac{\der L}{\der \dot{q}^{i}} \, , \\
\frac{\partial \tilde{L}}{\partial q^j}&=\frac{\partial L}{\partial q^j}+\left(y^a\frac{\partial X^i_a}{\partial q^j}+y^{\alpha}\frac{\partial X^i_\alpha}{\partial q^j}\right)\frac{\partial L}{\partial \dot{q}^i} \, .
\end{align*}
Now define the restricted Lagrangian $l: {\mathcal D} \rightarrow {\mathbb R}$ by $l:=\tilde{L}\big|_{\mathcal D}$, that is,  $l(q^i, y^a):=\tilde{L}(q^i, y^a, 0)$. It is interesting to note that 
\begin{equation*}\label{lagrangian}
\tilde{L}(q^i, y^a, y^{\alpha})=\frac{1}{2} g_{ab}y^ay^b+\frac{1}{2}g_{\alpha\beta}y^{\alpha}y^{\beta}-V(q) \, ,
\end{equation*}
where $g_{ab}:=g(X_a, X_b)$ and $g_{\alpha\beta}:=g(X_{\alpha}, X_{\beta})$, and thus
\[
l(q^i, y^a)=\frac{1}{2} g_{ab}y^ay^b-V(q) \, .
\]
We will make use of the fact that we can express the bracket $[X_a, X_b]$ in two ways using the different frames, concretely as 
\begin{align*}%\label{ers}
\left[X_a, X_b\right]&=\left(\frac{\partial X^j_b}{\partial q^i}X^i_a-\frac{\partial X^j_a}{\partial q^i}X^i_b\right)\frac{\partial}{\partial q^j}=[X_a, X_b]^j\frac{\partial}{\partial q^j} \, ,\\
\left[X_a, X_b \right]&= {\mathcal C}_{ab}^c X_c+ {\mathcal C}_{ab}^{\alpha} X_{\alpha}\; .
\end{align*}
Now, taking the restriction of Equation (\ref{equation1}) to ${\mathcal D}$, that is, using that $y^{\alpha}=0$, we get
\begin{align*}%\label{cuenta}
0&=
X^i_a\frac{d}{dt}\lp\frac{\der L}{\der\dot q^{i}}\rp-X^i_a\frac{\der L}{\der q^{i}}\\
&=
\frac{d}{dt}\lp\frac{\der L}{\der\dot q^{i}}X^i_a \rp   -     \frac{d X^i_a}{dt}\frac{\der L}{\der\dot q^{i}}
-
X^i_a\frac{\der L}{\der q^{i}}\\
&=
\frac{d}{dt}\lp\frac{\der l}{\der  y^a} \rp   + [X_a, X_b]^i y^b \frac{\der L}{\der\dot q^{i}}
-
X^i_a\frac{\der {l}}{\der q^{i}} \, .
\end{align*}
The middle term of the last equation is
\begin{align*}
 [X_a, X_b]^i y^b \frac{\der L}{\der\dot q^{i}}&=( {\mathcal C}_{ab}^c X^i_c+ {\mathcal C}_{ab}^{\alpha} X^i_{\alpha})y^b \frac{\der L}{\der\dot q^{i}}\\
 &={\mathcal C}_{ab}^c y^b\frac{\der l}{\der y^c}+{\mathcal C}_{ab}^{\alpha}y^b\frac{\der \tilde{L}}{\der  y^{\alpha}}\\
 &={\mathcal C}_{ab}^c y^b\frac{\der l}{\der y^c} \, ,
\end{align*}
since $\partial\tilde{L}/\partial y^{\alpha}=g_{\alpha\beta}y^{\beta}=0$ because $y^{\alpha}=0$. 
In conclusion we have that 
\begin{align*}%\label{cuenta-2}
0&=
X^i_a\frac{d}{dt}\lp\frac{\der L}{\der\dot q^{i}}\rp-X^i_a\frac{\der L}{\der q^{i}}\\
&=
\frac{d}{dt}\lp\frac{\der l}{\der y^a} \rp   + {\mathcal C}_{ab}^c y^b\frac{\der l}{\der  y^c}
-
X^i_a\frac{\der {l}}{\der q^{i}} \, .
\end{align*}
Therefore, the equations of motion of the nonholonomic system are rewritten in terms of the restricted Lagrangian $l$  as 
\begin{subequations}
\label{lagEq}
\begin{align}
&\displaystyle\frac{d}{dt}\lp\frac{\der l}{\der y^a} \rp   +   {\mathcal C}_{ab}^c y^b\frac{\der l}{\der  y^c}
-
X^i_a\frac{\der {l}}{\der q^{i}}=0\, , \label{lageq1}\\
&\dot{q}^i= X^i_a(q)y^a\, , \label{lageq2}
\end{align}
\end{subequations}
see for instance \cite{MR2492630,neimarkfufaev}.

\medskip

\hspace{0.015\linewidth}
\begin{minipage}{0.95\textwidth}
\begin{ex*}[Rolling disk, continued]
{	\rm \small 
We have the restricted Lagrangian 
\[
l(x_1, x_2, \theta, \phi, y^1, y^2)=\frac{1}{2}\left[(mr^2+J_{\phi})\,(y^1)^2 +J_{\theta}(y^2)^2\right]\, .   
\]
Now observe that
\begin{align*}
[X_1, X_2]&=r\sin\theta\frac{\partial}{\partial x_1}-r\cos\theta\frac{\partial}{\partial x_2}\\
&= mr\sin\theta\, X_3 - mr\cos\theta\, X_4 \, .
\end{align*}
Therefore in this simple example we have $C_{ab}^c=0$ for all $1\leq a, b, c\leq 2$.
The equations of motion \eqref{lagEq} for this nonholonomic system are
\begin{eqnarray*}
&& \dot{x}_1=ry^1\cos\theta  \, , \quad \dot{\theta}=y^2  \, , \quad \dot{y}^1=0 \, , \\
&& \dot{x}_2=ry^1\sin\theta  \, , \quad \dot{\phi}=y^1 \, , \quad \dot{y}^2=0  \, , 
\end{eqnarray*}
which are immediately explicitly integrated. \hfill $\Diamond$
}
\end{ex*}
\end{minipage}

\subsection{``Hamiltonian equations" for nonholonomic systems}
On the cotangent bundle $T^*Q$ the Lagrangian is replaced by the corresponding Hamiltonian $H$. We still assume a mechanical system, and let $(q^i,p_i)$, $1 \leq i \leq n$, give local canonical coordinates on $T^*Q$ through the Legendre transformation $\mathcal{F}L:TQ \rightarrow T^*Q$, i.e. 

\[
\mathcal{F}L:\, (q^i,\dot{q}^i)\, \longmapsto\, (q^i,p_i = \partial L/\partial \dot{q}^i)\, .
\]
Then $H$ is given locally by
\[
H(q^i,p_i) = \frac{1}{2}p_i\textbf{\textit{g}}^{ij}p_j + V(q)\, ,
\]
with $(\textbf{\textit{g}}^{ij})$ being the inverse matrix of $(\textbf{\textit{g}}_{ij})$.  

The Hamiltonian form of the nonholonomic equations \eqref{LdAeqs} is then

\begin{subequations}
\label{eq:HamEq}
\begin{align}
\begin{pmatrix}
\dot{q} \\
\dot{p}
\end{pmatrix}
&=
J
\left(
\begin{array}{c}
\frac{\partial H}{\partial q}(q,p) \\
\frac{\partial H}{\partial p}(q,p)
\end{array}
\right)
+
\lambda_\alpha
\begin{pmatrix}
0 \\
\mu^\alpha (q)
\end{pmatrix}\, ,  \label{eq:HamDyn}  \\
\mu^\alpha_i (q)\frac{\partial H}{\partial p_i}(q,p) &= \mu_{i}^{\alpha}\textbf{\textit{g}}^{ik}p_k = 0\, ,  \label{eq:HamConstraint} 
\end{align}
\end{subequations}
where $m+1 \leq \alpha \leq n$ and
$J=\left( 
\begin{array}{cc}
0_n&I_n\\
-I_n&0_n
\end{array}
\right)
$, see e.g. \cite{SM1994}.

\subsubsection{Equations of motion in adapted coordinates}
Now we will rewrite the restricted nonholonomic equations in a Hamilton-like way on $D^*$ (see \cite{Balseiro,hj-noholonomic,MR2492630}). More precisely, 
consider the Legendre transformation 
$
{\mathcal F}l:{\mathcal D}\rightarrow {\mathcal D}^* \, ,
$
locally given by 
\[
{\mathcal F}l:\, (q^i, y^a) \, \longmapsto\,  (q^i, \rho_a=\frac{\partial l}{\partial y^a})\, ,
\]
From the Legendre transformation we can define the Hamiltonian function $\mathcal{H}: {\mathcal D}^*\rightarrow {\mathbb R}$ 
which in local coordinates becomes
\begin{equation*}
\mathcal{H}(q^i, \rho_a)=\frac{1}{2}g^{ab}\rho_a\rho_b+V(q) \, . \label{hamilt}
\end{equation*}
Then upon changing coordinates in \eqref{lagEq} using the Legendre transform and $\mathcal{H}$, the equations of motion of a nonholonomic system are equivalently rewritten as 
\begin{subequations}
\label{H}
\begin{align}
\dot{q}^i&=X^i_b\frac{\partial \mathcal{H}}{\partial \rho_b}\; , \label{H1} \\
\dot{\rho}_a&=-{\mathcal C}_{ab}^c \rho_c\frac{\partial \mathcal{H}}{\partial \rho_b}-X^i_a\frac{\partial \mathcal{H}}{\partial q^i}\; . \label{H2}
\end{align}
\end{subequations}
If we define the skew-symmetric matrix
\begin{equation}\label{PiMat}
\Pi(q,\rho)=
\left(
\begin{array}{cc}
0 & X^i_b \\
-(X^j_a)^T & -C_{ab}^c \rho_c
\end{array}
\right)
\end{equation}
then the equations (\ref{H}) will be given by
\begin{equation}\label{dynamical}
\dot{\zeta}=\Pi(\zeta)\nabla \mathcal{H}(\zeta)\, ,
\end{equation}
where $\zeta=(q^i, \rho_a)$ are coordinates on ${\mathcal D}^*$. This skew gradient format will allow the use of discrete gradient methods, as we will see in the next section.

\begin{remark}
It is possible to give a more intrinsic definition of these objects, as in \cite{hj-noholonomic}.  Denote by  $\{\cdot, \cdot\}$ the canonical bracket of the cotangent bundle $T^*Q$. Define 
a bracket of functions $\{\cdot, \cdot\}_{{\mathcal D}^*}$ on $D^*$ by
\[
\{f, g\}_{{\mathcal D}^*} = \{f\circ i^*_{\mathcal D}, 
g\circ i^*_{\mathcal D} \}
\circ  P^*\, ,
\]
for $f, g\in C^{\infty}({\mathcal D}^*)$
 where $i^*_{\mathcal D}: T^*Q\rightarrow {\mathcal D}^*$ and $P^*: 
 {\mathcal D}^*\rightarrow T^*Q$ are the dual maps of the monomorphisms
$i_{\mathcal D}:  {\mathcal D}\rightarrow TQ$ and the projector $P: TQ\rightarrow {\mathcal D}$, respectively.
Then the bivector field $\Pi$ is given by
\[
\Pi(df, dg)=\{f, g\}_{{\mathcal D}^*} \, .
\]
This bracket does not in general satisfy the Jacobi identity, that is 
\[
\{f, \{g, h\}_{{\mathcal D}^*}\}_{{\mathcal D}^*}+\{g, \{h, f\}_{{\mathcal D}^*}\}_{{\mathcal D}^*}+\{h, \{f, g\}_{{\mathcal D}^*}\}_{{\mathcal D}^*}\not=0 \, .
\]
By using this bracket, Equation (\ref{dynamical}) will be more appropriately written as 
\[
\dot{f}=\{f, \mathcal{H} \}_{{\mathcal D}^*} \hbox{   for all   }  f\in C^{\infty}({\mathcal D}^*) \, .
\]
\end{remark}

\medskip

\hspace{0.015\linewidth}
\begin{minipage}{0.95\textwidth}
\begin{ex*}[Rolling disk, continued]
{	\rm \small
We have the Hamiltonian function
\[
\mathcal{H}(x_1, x_2, \theta, \phi, \rho_1, \rho_2)=\frac{1}{2}\left[\frac{\rho_1^2}{mr^2+J_{\phi}}\, +\frac{\rho_2^2}{J_{\theta}}\right]   \, ,
\]
where $\rho_1=(mr^2+J_{\phi}) y^1$ and $\rho_2=J_{\theta} y^2$. The  skew-symmetric matrix (almost-Poisson structure) \eqref{PiMat} is given by
\[
\Pi=\left(
\begin{array}{cccccc}
0&0&0&0&r\cos\theta&0\\
0&0&0&0&r\sin\theta&0\\
0&0&0&0&0&1\\
0&0&0&0&1&0\\
-r\cos\theta&-r\sin\theta&0&-1&0&0\\
0&0&-1&0&0&0
\end{array}
\right)\, ,
\]
and the equations of motion \eqref{dynamical} are
\[
\left(
\begin{array}{c}
\dot{x}_1\\
\dot{x}_2\\
\dot{\theta}\\
\dot{\phi}\\
\dot{\rho}_1\\
\dot{\rho}_2
\end{array}
\right)
=
\left(
\begin{array}{cccccc}
0&0&0&0&r\cos\theta&0\\
0&0&0&0&r\sin\theta&0\\
0&0&0&0&0&1\\
0&0&0&0&1&0\\
-r\cos\theta&-r\sin\theta&0&-1&0&0\\
0&0&-1&0&0&0
\end{array}
\right)
\left(
\begin{array}{c}
0\\
0\\
0\\
0\\
\displaystyle{\frac{\rho_1}{mr^2+J_{\phi}}}\\
\displaystyle{\frac{\rho_2}{J_{\theta}}}
\end{array}
\right) \, .
\]
\hfill $\Diamond$
}
\end{ex*}
\end{minipage}

\section{Energy-preserving integrators based on discrete gradients}\label{section3}

In the previous section, we reduced the study of the nonholonomic dynamics to a system of differential equations  
\[
\dot{\zeta}=\Pi(\zeta)\nabla \mathcal{H}(\zeta)\, 
\]
on ${\mathcal D}^*$. In this section we will assume that $Q$ is a real vector space, therefore ${\mathcal D}^{*}\cong \R^{N}$ where $n+m=N$. For a generalisation to the case of Lie groups and homogeneous manifolds see \cite{CeOw}.

Since  nonholonomic dynamics does not preserve the almost-Poisson structure $\Pi$ in general,  we will focus on the preservation of the energy using geometric integrators which preserve exactly this quantity. In particular, we will use discrete analogues of the gradient of the Hamiltonian function 
\cite{McQuRo}. 

\subsection{Discrete gradients} \label{dg-examples}

For ODEs in skew-gradient form, i.e. $\dot{x}=\Pi(x)\nabla H(x)$ with $x\in \mathbb{R}^N$ and $\Pi(x)$ a skew-symmetric matrix, it is immediate to check that $H$ is a first integral. Indeed 
\[
\dot{H}=\nabla H(x)^T\dot{x}=\nabla H(x)^T \Pi(x)\nabla H(x)=0\, ,
\]
due to the skew-symmetry of $\Pi$. Using discretizations of the gradient $\nabla H(x)$ it is possible to define a class of integrators which preserve the first integral $H$ exactly.

\begin{definition}\label{def31}{\cite{GONZ}}
Let $H:\mathbb{R}^N\longrightarrow \mathbb{R}$ be a differentiable function. Then $\bar{\nabla}H:\mathbb{R}^{2N}\longrightarrow \mathbb{R}^N$ is a discrete gradient of $H$ if it is continuous and satisfies
\begin{subequations}
\label{discGrad}
\begin{align}
\bar{\nabla}H(x,x')^T(x'-x)&=H(x')-H(x)\, , \quad \, \mbox{ for all } x,x' \in\mathbb{R}^N  \, ,\label{discGradEn} \\
\bar{\nabla}H(x,x)&=\nabla H(x)\, , \quad \quad \quad \quad \mbox{ for all } x \in\mathbb{R}^N  \, . \label{discGradCons}
\end{align}
\end{subequations}
\end{definition}

Some well-known examples of discrete gradients are:
\begin{itemize}
\item The mean value (or averaged vector field) discrete gradient introduced in \cite{HLvL83} and given by
\begin{equation}
\label{AVF}
\bar{\nabla}_{1}H(x,x'):=\int_0^1 \nabla H ((1-\xi)x+\xi x')d\xi \, , \quad \mbox{ for } x'\not=x \, .
\end{equation}

\item The midpoint (or Gonzalez) discrete gradient, introduced in \cite{GONZ} and given by
\begin{align}
\bar{\nabla}_{2}H(x,x')&:=\nabla H\left( \frac{1}{2}(x'+x)\right)+\frac{H(x')-H(x)-\nabla H\left( \frac{1}{2}(x'+x)\right)^T(x'-x)}{|x'-x|^2}(x'-x) \, , \label{gonzales}\\
& \mbox{ for } x'\not=x \, . \nonumber
\end{align}
\item The coordinate increment discrete gradient, introduced in \cite{ITOH}, with each component given by
\begin{equation}
\label{itoAbe}
\bar{\nabla}_{3}H(x,x')_i=\frac{H(x'_1,\ldots,x'_i,x_{i+1},\ldots,x_n)-H(x'_1,\ldots,x'_{i-1},x_{i},\ldots,x_n)}{x'_i-x_i}\, , \quad 1\leq i \leq N\, ,
\end{equation}
when $x_i'\not=x_i$, and $\bar{\nabla}_{3}H(x,x')_i=\frac{\partial H}{\partial x_i}(x'_1,\ldots,x'_{i-1},x'_i=x_{i},x_{i+1},\ldots,x_n)$ otherwise.
\end{itemize}

It can be easily checked that these are indeed discrete gradients, see \cite{GONZ}, \cite{ITOH} and \cite{McQuRo}.

\subsection{Integrators based on discrete gradients}
Once a discrete gradient $\bar{\nabla}H$ has been chosen, it is straightforward to define an energy-preserving integrator by
\begin{equation}
\frac{x'-x}{h}=\tilde{\Pi}(x,x',h)\bar{\nabla}H(x,x') \, ,  \label{E-int}
\end{equation}
where $\tilde{\Pi}$ is a differentiable skew-symmetric matrix approximating $\Pi$, that is, it satisfies $\tilde{\Pi}(x,x,0)=\Pi(x)$. As in the continuous case, it is immediate to check that $H$ is exactly preserved, since
$$
H(x')-H(x)=\bar{\nabla}H(x,x')^T(x'-x)=h \bar{\nabla}H(x,x')^T \tilde{\Pi}(x,x',h)\bar{\nabla}H(x,x')=0 \, . 
$$
If we further wish to get a second order method then it is sufficient to choose $\tilde{\Pi}$ such that $\tilde{\Pi}(x,x',h)=\tilde{\Pi}(x',x,-h)$, and a differentiable discrete gradient such that $\bar{\nabla}H(x,x')=\bar{\nabla}H(x',x)$. This guarantees that the integration method (\ref{E-int}) is time-symmetric and therefore second order accurate, see \cite{hairer}. For instance it is enough to choose $\tilde{\Pi}(x,x',h)=\Pi\left( \frac{x+x'}{2} \right)$ and take the mean value discrete gradient or the midpoint discrete gradient. 
Higher order energy-preserving methods, which generalize the mean value discrete gradient \eqref{AVF}, can be obtained by collocation methods as in \cite{Hairer2011}.

\begin{remark}
\label{rem:IM}
If the Hamiltonian is quadratic then
$$
\bar{\nabla}_{1}H(x,x')=\bar{\nabla}_{2}H(x,x')=\nabla H\left( \frac{1}{2}(x'+x)\right) \, ,
$$
that is, the mean value discrete gradient \eqref{AVF} and the Gonzalez discrete gradient \eqref{gonzales} coincide with the continuous gradient evaluated at the midpoint. Then if we choose $\tilde{\Pi}(x,x',h)=\Pi\left( \frac{x+x'}{2} \right)$ the method (\ref{E-int}) reduces to the implicit midpoint rule.
If the Hamiltonian is of the form $H(x)=\sum_{j=i}^{N}a_jx_j^2$, then $(\nabla H)_j=2a_jx_j$ and
$$
\bar{\nabla}_{1}H(x,x')_i=\bar{\nabla}_{2}H(x,x')_i=\bar{\nabla}_{3}H(x,x')_i=a_i(x'_{i}+x_i) \, , \quad 1\leq i \leq N\, ,
$$
that is, all three discrete gradients introduced above coincide with $\nabla H\left( \frac{1}{2}(x'+x)\right)$.
\end{remark}

\begin{remark}{\bf Preservation of the nonholonomic constraints.}
Going back to the case of nonholonomic systems, we can now apply an energy-preserving method (\ref{E-int}) to Equation (\ref{dynamical}). Notice that if we take the approximation $\tilde{\Pi}(\zeta,\zeta',h)$ to be $\Pi({\bar{\zeta}})$ for some $\bar{\zeta}(\zeta,\zeta')\in D^*$, and let $\partial \bar{\mathcal{H}}/\partial \rho_a$ be the discrete gradient component that approximates $\partial \mathcal{H}/\partial \rho_a$, then a discrete gradient method \eqref{E-int} gives
$$
\frac{q'^j-q^j}{h}=X^j_a(\bar{q})\frac{\partial \bar{\mathcal{H}}}{\partial \rho_a}(\zeta,\zeta') \, ,
$$
where $\bar{q}=\pi_Q(\bar{\zeta})$.
When applying $\mu^\alpha\in D^\circ$ we obtain
$$
\mu_j^\alpha \left( \frac{q'^j-q^j}{h} \right)=\mu_j^\alpha X_a^j(\bar{q})\frac{\partial \bar{\mathcal{H}}}{\partial \rho_a}(\zeta,\zeta')=0 \, ,
$$
since $\mu_j^\alpha X_a^j=0$ for all $1 \leq a \leq m$, $m+1\leq \alpha \leq n$. All the nonholonomic constraints are thus preserved by the method. 
\end{remark}

\medskip

\hspace{0.015\linewidth}
\begin{minipage}{0.95\textwidth}
\begin{ex*}[Rolling disk, continued]
{	\rm \small
Using any of the three discrete gradients introduced in Section \ref{dg-examples} and a midpoint approximation of $\Pi$, we get the following energy preserving integrator, which is precisely the implicit midpoint rule:
\begin{align*}
x_1'&=x_1+hr\cos\left(\frac{\theta+\theta'}{2}\right)\frac{\rho_1}{mr^2+J_{\phi}} \, , \quad \theta'=\theta+h\frac{\rho_2}{J_{\theta}} \, , \quad \rho'_1=\rho_1 \, , \\
x_2'&=x_2+hr\sin\left(\frac{\theta+\theta'}{2}\right)\frac{\rho_1}{mr^2+J_{\phi}} \, , \quad \phi'=\phi+h\frac{\rho_1}{mr^2+J_{\phi}} \, , \quad \rho'_2=\rho_2 \, .  
\end{align*}
Observe that, as a consequence, we deduce the preservation of the constraints
\[
\frac{x_1'-x_1}{h}-r\left( \frac{\phi'-\phi}{h}\right)\cos\left( \frac{\theta+\theta'}{2}\right)=0\; , \quad \frac{x_2'-x_2}{h}-r \left(\frac{\phi'-\phi}{h}\right)\sin\left( \frac{\theta+\theta'}{2}\right)=0 \, , 
\]
which are discretizations of the nonholonomic constraints \eqref{diskConst}.
\hfill $\Diamond$
}
\end{ex*}
\end{minipage}

\subsection{Integrators on $T^*Q$} \label{T*Q}
The equations of motion in adapted coordinates for a given nonholonomic system are usually not known initially. A potential obstacle in applying a discrete gradient method to the equations in adapted coordinates \eqref{dynamical} is therefore that the user must analytically derive these equations. In this section we formulate the proposed schemes directly on the Hamiltonian equations of motion in canonical coordinates \eqref{eq:HamEq}, and achieve preservation of energy and the nonholonomic constraints without explicitly deriving and using the formulation in adapted coordinates.

As a first attempt at an energy preserving method, we can define a numerical integrator for \eqref{eq:HamEq} directly on $T^*Q$ by
\begin{equation} \label{simpDirect}
\frac{z'-z}{h}
=
J
\bar{\nabla}H(z,z')
+
\lambda_\alpha
\begin{pmatrix}
0 \\
\mu^\alpha (\bar{q})
\end{pmatrix}
\, , \quad 
\bar{\nabla}H(z,z')^T \begin{pmatrix}
0 \\
\mu^\alpha (\bar{q})
\end{pmatrix}=0
\, \, ,
\end{equation}
where $\bar{\nabla}H$ is a discrete gradient, $z=(q, p)$ and $z'=(q',p')$. Notice that this method is energy-preserving, since
\begin{align*}
H(z')-H(z)&=\bar{\nabla}H(z,z')^T(z'-z) \\
&=h\bar{\nabla}H(z,z')^T J\bar{\nabla}H(z,z')+h\lambda_\alpha\bar{\nabla}H(z,z')^T \begin{pmatrix}
0 \\
\mu^\alpha (\bar{q})
\end{pmatrix}=0 \, .
\end{align*}
However the constraints \eqref{eq:HamConstraint} will in general only be approximately satisfied at the solution points by such a method.

To achieve exact preservation of both the energy and the nonholonomic constraints (\ref{eq:HamConstraint}), we utilize the restricted equations \eqref{dynamical} on $\mathcal{D}^{*}$.

The method to step from $(q,p)$ to $(q',p')$ can be summarized as:
\begin{enumerate}
\item Change coordinates from $(q,p)$ to $(q,\rho)$.
\item Step from $(q,\rho)$ to $(q',\rho')$ using a discrete gradient method \eqref{E-int} applied to \eqref{dynamical}.
\item Change coordinates from $(q',\rho')$ to $(q',p')$.
\end{enumerate}

For step (i) and (iii) we make use of the following relations between the coordinates $\rho_b$ on $\mathcal{D}^*$ and $p_i$ on $\mathcal{F}L(\mathcal{D})$
\begin{subequations}
\label{ptransform}
\begin{align}
p_i &= \textbf{\textit{g}}_{ij}X_a^jg^{ab}\rho_b \, , \quad  1\leq i \leq n \, , \label{ptoptilde}\\ 
\rho_b &= X_b^i p_i \, , \quad 1\leq b \leq m \, . \label{ptildetop}
\end{align} 
\end{subequations}
The challenge is performing step (ii) without any explicit knowledge of the equations in adapted coordinates. Specifically we need to evaluate $\nabla \mathcal{H}$ and $\Pi$ in \eqref{dynamical} for any $\zeta \in \mathcal{D}^*$. Let us therefore rewrite these expressions in a suitable format. Here $p$ is considered a dependent variable of $(q,\rho)$ through \eqref{ptoptilde}.

First observe that the skew symmetric matrix $-C_{ab}^c\rho_c$ in \eqref{PiMat} may be written as
\begin{equation}
-C_{ab}^c\rho_c = -C_{ab}^cX_c^jp_j  = -[X_a,X_b]^jp_j = \left(\frac{\partial X_a^j}{\partial q^i}X_b^i-\frac{\partial X_b^j}{\partial q^i}X_a^i\right)p_j \, . 
\label{skewMat}
\end{equation}
Second we can write the partial derivatives of $\mathcal{H}$ as follows
\begin{subequations}
\label{eq:Htgrad}
\begin{align}
\frac{\partial \mathcal{H}}{\partial q^i} &= \frac{1}{2}p_j\frac{\partial \textbf{\textit{g}}^{jk}}{\partial q^i}p_k+ \frac{\partial p_j}{\partial q^i} \textbf{\textit{g}}^{jk}p_k+\frac{\partial V}{\partial q^i} = \frac{1}{2}p_j\frac{\partial \textbf{\textit{g}}^{jk}}{\partial q^i}p_k+ \frac{\partial p_j}{\partial q^i}X_a^jg^{ab}\rho_b+\frac{\partial V}{\partial q^i} \nonumber \\
&= \frac{1}{2}p_j\frac{\partial \textbf{\textit{g}}^{jk}}{\partial q^i}p_k- p_j\frac{\partial X_a^j}{\partial q^i}g^{ab}\rho_b+\frac{\partial V}{\partial q^i} \, , \label{eq:Htgradq}  \\
\frac{\partial \mathcal{H}}{\partial \rho_a} &=  g^{ab}\rho_{b} = g^{ab}X_b^ip_i \, . \label{eq:Htgradp}
\end{align}
\end{subequations}

Expressing \eqref{dynamical} using \eqref{skewMat} and \eqref{eq:Htgrad}, the remaining issue is that we don't have explicit knowledge of a basis $X_a(q)$ for the distribution $\mathcal{D}$ or of the partial derivatives $\partial X_a(q)/\partial q^i$. For an arbitrary point $q$ we generate $X_a(q)$ by computing the QR-factorization of the constraint matrix $(\mu_{i}^{\alpha}(q))$ using Householder reflections, see e.g. \cite{golubvanloan}. The last $m$ columns of the $Q$ matrix can then be taken as $X_a(q)$, $1 \leq a \leq m$. Householder reflections were chosen because they are numerically stable and efficient. 

We now make the assumption that the partial derivatives of $(\mu_{i}^{\alpha}(q))$ are either known or easily derived, which is usually the case. Then $\partial X_a(q)/\partial q^i$ can be calculated by augmenting the QR-factorization algorithm with corresponding steps for the partial derivatives. The procedure is specified in Algorithm \ref{algQRdiff} for the matrix $A(q) := (\mu_{i}^{\alpha}(q))$. 

\begin{algorithm}
\small
\begin{multicols}{2}
\begin{algorithmic}[0]
\Procedure{QRdiff}{$A,\partial A,n,m,s$}
\State $Q^{(0)} \gets I_n$
\State $R^{(0)} \gets A$
\For{$i=1,2,\hdots,n$}\
\State $\partial_i Q^{(0)} \gets 0_n$
\State $\partial_i R^{(0)} \gets \partial_i A$
\EndFor
\Statex
\For{$k=0,1,\hdots,n-m-1$}\
\For{$j=1,2,\hdots,n$}\
\If{$j\leq k$}
\State $\tilde{w}^{(k)}_j \gets 0$
\ElsIf{$j = k+1$}
\State $\tilde{w}^{(k)}_j \gets R^{(k)}_{jk} + s_k\|R^{(k)}_{k:n,k}\|_2$
\Else
\State $\tilde{w}^{(k)}_j \gets R^{(k)}_{jk}$
\EndIf
\EndFor
\State $\|\tilde{w}^{(k)}\|_2 \gets \sqrt{2\left(\|R^{(k)}_{k:n,k}\|^2_2+s_k R^{(k)}_{kk}\|R^{(k)}_{k:n,k}\|_2\right)} $
\State $w^{(k)} \gets \frac{\tilde{w}^{(k)}}{\|\tilde{w}^{(k)}\|_2}$
\State $u^{(k)} \gets {w^{(k)}}^TR^{(k)}$
\State $R^{(k+1)} \gets R^{(k)}-2w^{(k)}u^{(k)} $
\State $Q^{(k+1)} \gets Q^{(k)}-2Q^{(k)}\left(w^{(k)}{w^{(k)}}^T \right)$
  \For{$i=1,2,\hdots,n$}\
  \State $ \partial_i \|R^{(k)}_{k:n,k}\|_2 \gets \frac{\left(\partial_i R^{(k)}_{k:n,k}\right)^T R^{(k)}_{k:n,k}}{\|R^{(k)}_{k:n,k}\|_2} $
\For{$j=1,2,\hdots,n$}\
\If{$j\leq k$}
	\State $\partial_i\tilde{w}^{(k)}_j \gets 0$
\ElsIf{$j = k+1$}
	\State $\partial_i\tilde{w}^{(k)}_j \gets \partial_i R^{(k)}_{jk} + s_k\partial_i \|R^{(k)}_{k:n,k}\|_2$
\Else
	\State $\partial_i\tilde{w}^{(k)}_j \gets \partial_i R^{(k)}_{jk}$
\EndIf
\EndFor
 \State $\partial_i\|\tilde{w}^{(k)}\|_2 \gets \frac{\left(\partial_i \tilde{w}^{(k)}\right)^T \tilde{w}^{(k)}}{\|\tilde{w}^{(k)}\|_2} $
 \State $\partial_i w^{(k)} \gets \frac{\partial_i \tilde{w}^{(k)}\|\tilde{w}^{(k)}\|_2-\tilde{w}^{(k)} \partial_i \|\tilde{w}^{(k)}\|_2}{\|\tilde{w}^{(k)}\|^2_2} $
 \State $\partial_i u^{(k)} \gets {\partial_i w^{(k)}}^TR^{(k)}+{w^{(k)}}^T\partial_i R^{(k)}$
 \State $\partial_i R^{(k+1)} \gets \partial_i R^{(k)}-$
 \State $\quad 2\left(\partial_i w^{(k)}u^{(k)}+w^{(k)}\partial_i u^{(k)}\right)$
 \State $\partial_i Q^{(k+1)} \gets \partial_i Q^{(k)}-2\partial_i Q^{(k)}\left(w^{(k)}{w^{(k)}}^T\right)+$
 \State $\quad 2Q^{(k)}\left(\partial_i w^{(k)}{w^{(k)}}^T + w^{(k)}\partial_i{w^{(k)}}^T \right)$
\EndFor
\EndFor
\State \Return{$Q^{(n-m)},\partial Q^{(n-m)}$}
\EndProcedure
\end{algorithmic}
\end{multicols}
\captionof{algorithm}{(QR factorization procedure with differentiation using Householder reflections). Computes the QR factorization of a differentiable matrix $A(q)\in \mathbb{R}^{n,n-m}$ as well as all partial derivatives $\partial Q/\partial q^i$, $1 \leq i \leq n$, for any $q\in \mathbb{R}^n$. Let $\partial_i := \partial / \partial q^i$ while $B_{i:j,k:l}$ with $i\leq j$ and $k\leq l$ denotes the submatrix containing rows $i$ to $j$ and columns $k$ to $l$ of a matrix $B$, with the shorthand $i := i:i$. $\partial B$ denotes the tensor containing all partial derivatives of $B$ at $q$. $s\in \{+1,-1\}^{n-m}$ is a vector of sign choices.}
\label{algQRdiff}
\end{algorithm} 

To ensure we are sampling the same basis vector fields at different points in a given step when using Algorithm \ref{algQRdiff}, it is sufficient to make sure the vector of sign choices $s$ remains fixed for all factorizations in a given integration step. Because we only suppose knowledge of the full system \eqref{eq:HamEq}, we transform back to canonical coordinates $(q,p)$ after each step.

In theory this implementation can be combined with any discrete gradient method. However, since it is desirable to minimize the number of QR-factorizations per time step, this approach is best suited when used together with the Gonzalez discrete gradient and a midpoint approximation of $\Pi$. We shall refer to this specific method later as GONZALEZ-R.
\begin{remark}{\bf Computational cost.}
For the initial direct method \eqref{simpDirect} it is necessary to evaluate the Lagrange multipliers. Moreover it is necessary to implement the constraint equations in each step of the algorithm. Applying a discrete gradient method \eqref{E-int} directly to the reduced system \eqref{dynamical} simplifies the computational cost. This is so because the constraints are preserved automatically, and it is not necessary to compute the Lagrange multipliers as additional variables. Specifically with the method \eqref{simpDirect} it is necessary to solve $3n-m$ variables while using \eqref{E-int} on \eqref{dynamical} it is only necessary to compute $n+m$ variables.

Integrating the full system using the equations in adapted coordinates and the QR-factorization approach, we avoid the problem with Lagrange multipliers, but still see a rise in computational cost due to the necessity of moving between coordinate systems, and the general added cost in evaluating $\nabla \mathcal{H}$ and $\Pi$.

Thus the trade-off in not requiring knowledge of the reduced system is an increase in computational cost. It is therefore generally more efficient to analytically derive \eqref{dynamical} and apply a discrete gradient method \eqref{E-int}.
\end{remark}
\begin{remark}{\bf Implementation using finite differences.} 
It is remarked in \cite{perlmutter06} that the condition \eqref{discGradCons} for discrete gradients is only required to ensure consistency. Suppose this condition is relaxed slightly to
\[
\bar{\nabla} H(x,x,h) = \nabla H(x) + \mathcal{O}(h^r) \, ,
\]
where $r$ should at least match the order of the method. This is sufficient for the consistency of an integrator \eqref{E-int}, and indeed for the method to have order $r$. We can use this to avoid having to evaluate $\partial X_a(q)/\partial q^i$ at the midpoint in GONZALEZ-R by replacing it with an appropriate finite difference approximation, e.g. the central difference approximation
\[
\bar{\frac{\partial X_a^j}{\partial q^i}}(q) := \frac{X_a^j(q+he_i)-X_a^j(q-he_i)}{2h} = \frac{\partial X_a^j}{\partial q^i}(q) + \mathcal{O}(h^2) \, ,
\]
where $e_i$ is the canonical unit vector $i$. The resulting method retains second order and still preserves energy and the nonholonomic constraints.
\end{remark}

\section{Examples and numerical experiments}\label{section4}
In this section we apply discrete gradient methods to some illustrative examples of nonholonomic systems.
In the first three examples we will derive Equations \eqref{dynamical} analytically. In the last one we will compare the strategies proposed in Section \ref{T*Q}.

\subsection{A fully chaotic nonholonomic system}\label{section41}

In \cite{perlmutter06} the authors remark that the key geometric properties for nonholonomic dynamics are not known for general nonintegrable systems. To compare integration methods for such systems, they focus on energy preservation, looking at the following chaotic quartic mechanical system on $Q = \mathbb{R}^{2n+1}$ with coordinates $q = (q^1,q^2,\hdots,q^{2n+1})^T:=(x,w_1,\hdots,w_n,z_1,\hdots,z_n)^T$, which is defined by the Lagrangian
\begin{equation}
\label{Lag}
L(q,\dot{q}) = K(q,\dot{q})-V(q), \ \ \mbox{where} \ \  K = \frac{1}{2}\|\dot{q}\|_2^2, \ \ \ \ V = \frac{1}{2}\left(\|q\|_2^2 + z_1^2z_2^2 + \sum_{i=1}^{n}w_i^2z_i^2\right) \, ,
\end{equation}
with the single nonholonomic linear velocity constraint
\begin{equation}
\label{constraint}
\dot{x} + \sum_{i=1}^{n} w_i\dot{z}_i = 0 \, .
\end{equation}
This system is reversible and preserves energy, i.e. $\dot{H} = 0$. 

\subsubsection{Lagrangian side}
To derive the equations in adapted coordinates, first note that this Lagrangian is of mechanical type. We can write the kinetic energy $K$ as 
\[
K(q,\dot{q}) = \frac{1}{2}g_{q}(\dot{q},\dot{q}) \, ,
\] with the canonical Riemannian metric
\[
g = \sum_{i=1}^{2n+1} dq^i \otimes dq^i \, ,
\]
which does not depend on $q$.

The distribution $\mathcal{D}$ and its orthogonal complement $\mathcal{D}^{\bot,g}$ are given by, respectively, the span of $2n$ and $1$ independent vector fields:
\begin{align*}
\mathcal{D} &= \textnormal{span}\left\{X_{i} := \frac{\partial}{\partial w_i}, \, X_{n+i} := w_i\frac{\partial}{\partial x}-\frac{\partial}{\partial z_i}, \, i=1,\hdots,n \right\}\, , \\
\mathcal{D}^{\bot,g} &= \textnormal{span}\left\{X_{2n+1} := \frac{\partial}{\partial x} + \sum_{i=1}^{n}w_i \frac{\partial}{\partial z_i} \right\} \, . \\
\end{align*}
The adapted basis $\{X_1,X_2,\hdots,X_{2n+1}\}$, induces new coordinates $(q^{i},y^{a},y^{2n+1})$, $1 \leq i \leq 2n+1$, $1 \leq a \leq 2n$ on $TQ$ for which the nonholonomic constraint reduces to $y^{2n+1} = 0$. The restricted Lagrangian $l: \mathcal{D} \rightarrow \mathbb{R}$ for this system is 
\[
l(q^i,y^a) = \frac{1}{2}\left(g_{ab}y^ay^b-V(q)\right) \, ,
\]
where 
\[
(g_{ab}) = \begin{pmatrix}
I_n& 0_n\\
0_n&I_n+ww^T\\
\end{pmatrix}  \quad \mbox{ and } \quad w = (w_1,w_2,\hdots,w_n)^T.
\]

\subsubsection{Hamiltonian side}
Moving to the Hamiltonian side we replace the velocities in adapted coordinates $y^a$ with the momenta
\[
\rho_{a} = \frac{\partial l}{\partial y^{a}} = g_{ab}y^b, \quad a=1,\hdots,2n \, .
\]
The restricted Hamiltonian $\mathcal{H}: \mathcal{D}^* \rightarrow \mathbb{R}$ is then
\[
\mathcal{H}(q^i,\rho_{a}) = \frac{1}{2} g^{ab}\rho_a\rho_b+V(q) \, , \\
\]
where $(g^{ab})$ is the inverse of $(g_{ab})$, i.e.
\[
(g^{ab}) = \begin{pmatrix}
I_n& 0_n\\
0_n&I_n-\frac{ww^T}{1+w^Tw}\\
\end{pmatrix} \, .
\]
The equations of motion on the Hamiltonian side are then given by \eqref{dynamical} where
\begin{equation}
\label{fcseq}
\Pi(q,\rho) = \begin{pmatrix}
& & & 0_{1\times n} & w^T \\
& 0_{2n+1} & &I_n & 0_n \\
& & & 0_n & -I_n \\
0_{n\times 1} & -I_n & 0_n & 0_n & -\kappa I_n \\
-w & 0_n & I_n & \kappa I_n & 0_n
\end{pmatrix}, \quad 
\nabla \mathcal{H}(q,\rho) = \begin{pmatrix}
x \\
w_1 + w_1z_1^2 -\kappa\rho_{n+1} + \kappa^2w_1 \\
\vdots \\
w_n + w_nz_n^2 -\kappa \rho_{n+n} + \kappa^2w_n \\
z_1 + y_1^2z_1 +z_1z_2^2 \\
z_2 + y_2^2z_2 +z_1^2z_2 \\
z_3 + y_3^2z_3  \\
\vdots \\
z_n + y_n^2z_n  \\
g^{ab}\rho_b
\end{pmatrix} \, , 
\end{equation}
$\kappa := \frac{w^T\eta}{1+w^Tw}$ and $\eta := (\rho_{n+1},\hdots,\rho_{2n})^T$.

\subsubsection{Numerical experiments}
We follow the approach in \cite{perlmutter06}, and integrate this system, with $n=3$, from a random initial state with energy $H=3.06$. We here compare five different methods. The first two are variational integrators based on the discrete Lagrange d'Alembert (DLA) for the full system \eqref{Lag} and \eqref{constraint}. The semi implicit reversible DLA variational integrator proposed in \cite{perlmutter06} (SI-DLA), and the implicit reversible DLA variational integrator based on a midpoint discrete Lagrangian (I-DLA) which is also described in \cite{perlmutter06} among others. The third method is the 2-stage Lobatto IIIA-B-C-C*-D SPARK method described in \cite{jay03} for index 2 DAEs (SPARK), which again discretize the equations of motion of the full system. For the last two methods, we integrate the reduced system, \eqref{dynamical} with \eqref{fcseq}, using a discrete gradient method \eqref{E-int}, with two different discrete gradients: The averaged vector field discrete gradient (AVF) \eqref{AVF}, and the Gonzalez discrete gradient (GONZALEZ) \eqref{gonzales}. 

As seen in Figure \ref{fig:FCS}, while all five methods are known to be second order accurate and respect the constraint, only the discrete gradient methods conserve the energy up to round off error. In \cite{perlmutter06} it is shown that the energy error for SI-DLA closely follows a random walk with the variance $\sigma^2 = 10^{-4}h^4t$. In Figure \ref{fig:FCS} we also show that I-DLA and SPARK behaves similarly, with all comparison methods exhibiting similar linear time growth. As expected, since the discrete gradient methods have no energy error, they also have zero variance up to round off error.

\begin{minipage}[t]{0.45\textwidth}
\centering
\includegraphics[width=\textwidth]{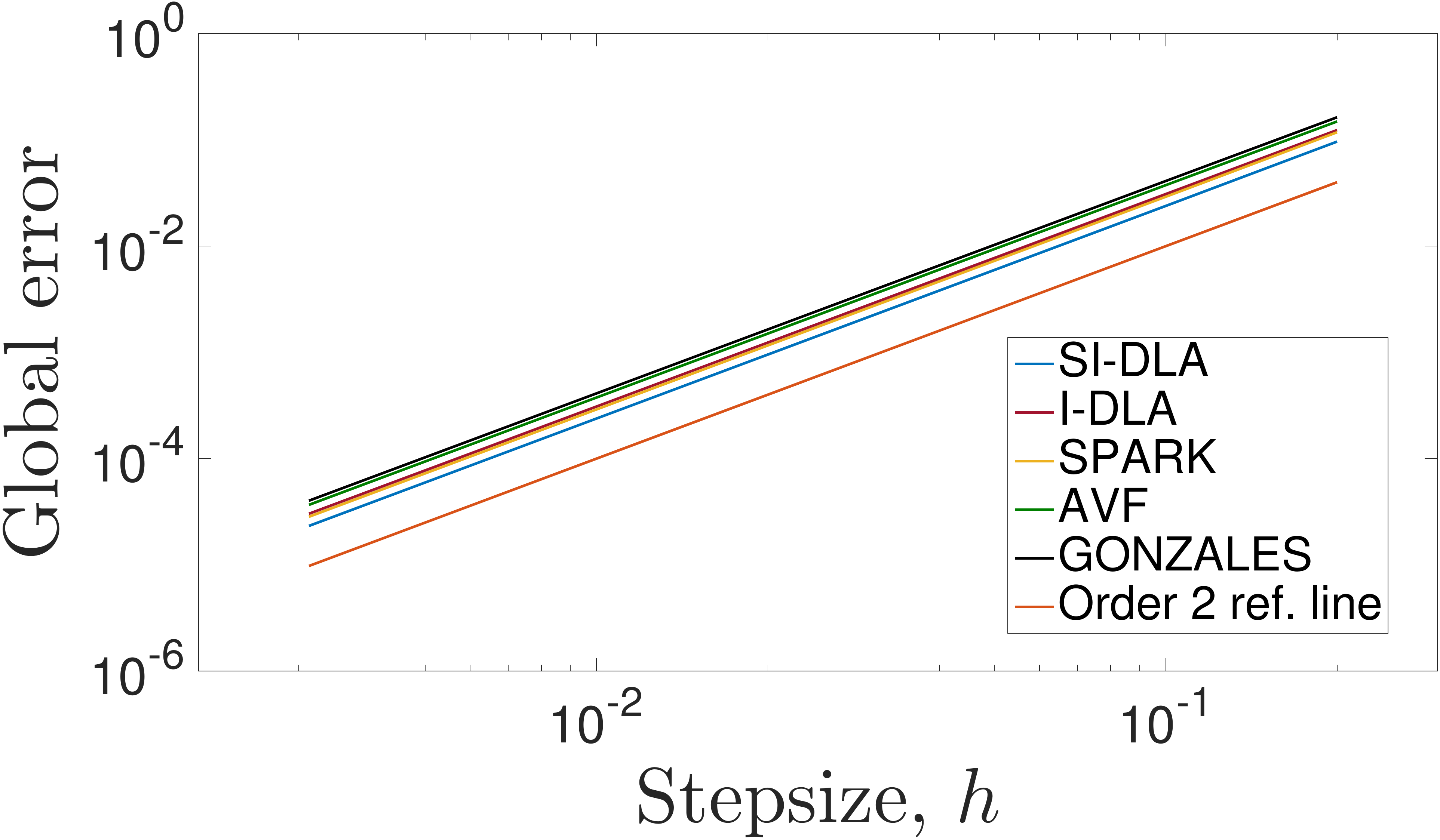}
\end{minipage}
\hspace{0.05\linewidth}
\begin{minipage}[t]{0.45\textwidth}
\centering
\includegraphics[width=\textwidth]{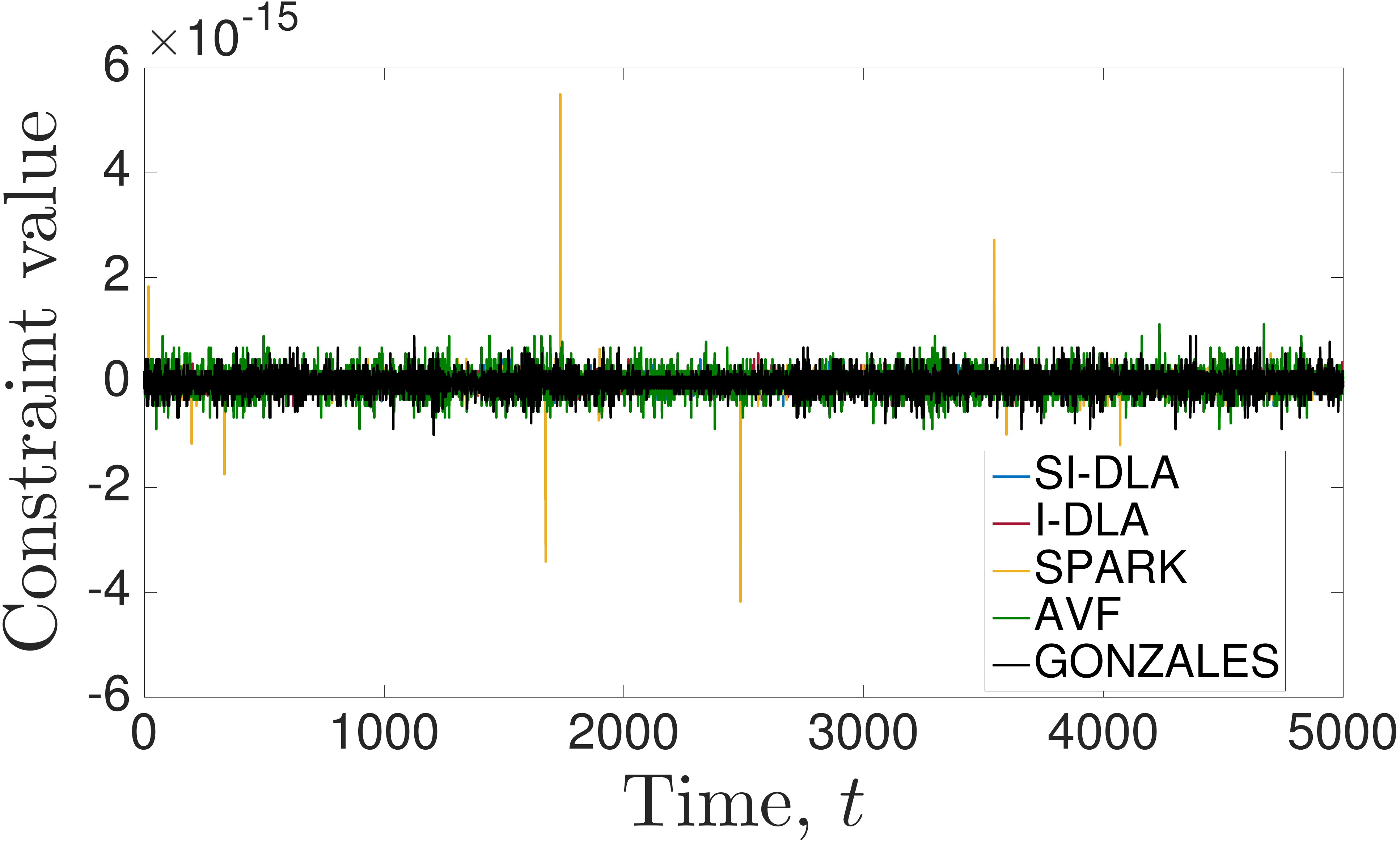}

\end{minipage}

\begin{minipage}[t]{0.45\textwidth}
\centering
   \includegraphics[width=\textwidth]{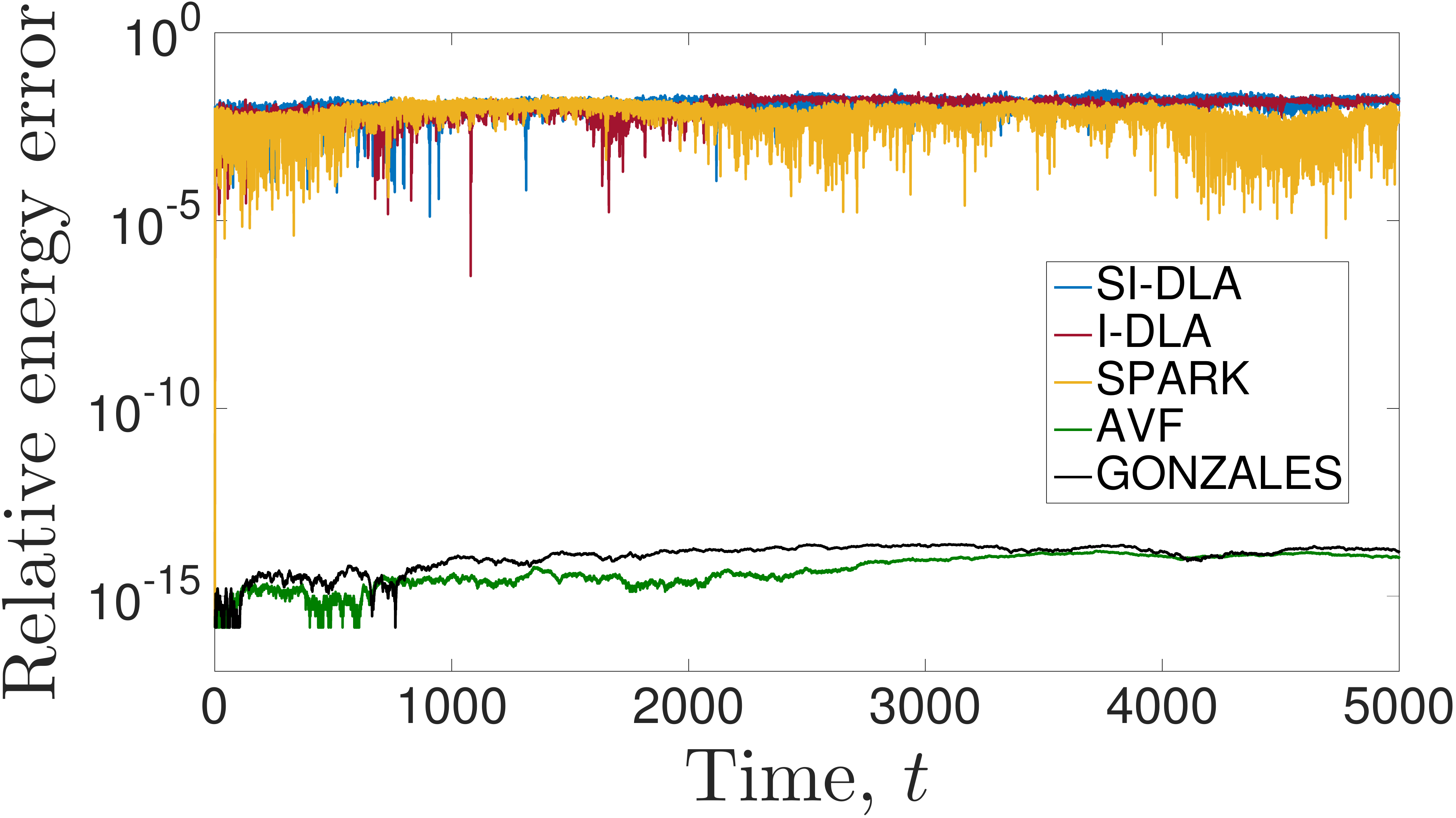}

\end{minipage}
\hspace{0.05\linewidth}
\begin{minipage}[t]{0.45\textwidth}
\centering
   \includegraphics[width=\textwidth]{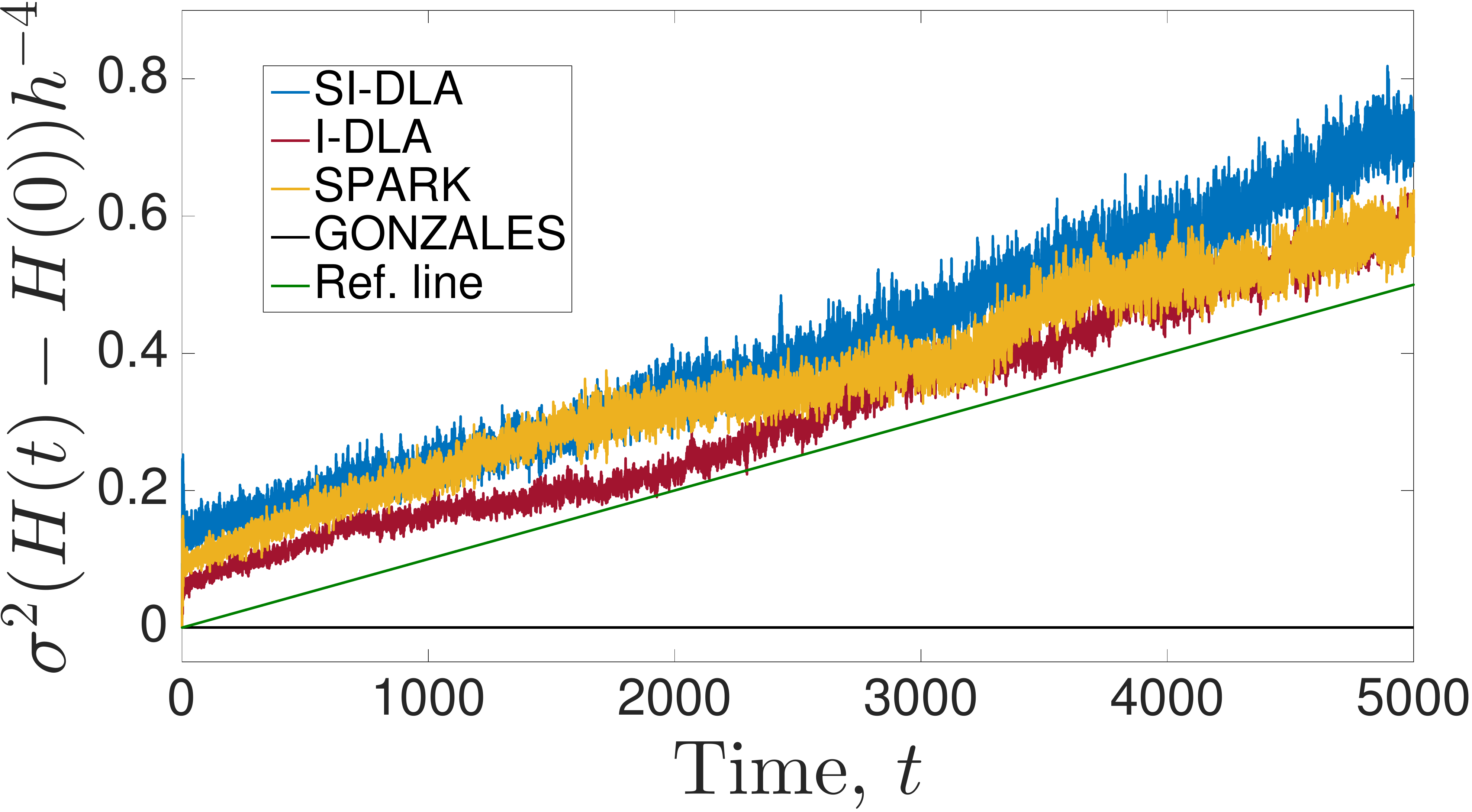}

\end{minipage}
\begingroup
 \captionof{figure}{\label{fig:FCS} {\footnotesize Comparison of the different methods for the fully chaotic system \eqref{Lag}-\eqref{constraint}. {\bf Top left}: Order plot, integrating up to $t=10$. All methods are seen to be second order. {\bf Top right}: Value of the left hand side of the constraint expression \eqref{constraint} for a sample trajectory with random initial conditions. The methods all respect the constraint up to machine precision. {\bf Bottom left}: Relative energy error, i.e. $|H(t)-H(0)/H(0)|$, for the same trajectory. Only the discrete gradient methods conserve the energy up to machine precision. {\bf Bottom right}: The variance of the energy errors $\sigma^2(H(t)-H(0))$ for 200 different initial conditions scaled by their expected $h^4$ dependence on the time step. The reference line $10^{-4}t$ is included for comparison. All comparison methods exhibit similar linear time growth in accordance with the reference line, while the discrete gradient method GONZALEZ has zero variance up to machine precision as expected. AVF is not shown since it was indistinguishable from GONZALEZ.}}
 \endgroup

\subsection{The Chaplygin sleigh}\label{section42}

In this example we will see that the transformation of the systems into adapted coordinates can give rise to some additional numerical advantages apart from the possibility of achieving energy preservation.

The Chaplygin sleigh is a rigid body moving on a horizontal plane with three contact points, two of which slide freely without friction. The third one is a knife edge, which imposes the nonholonomic constraint of no motion perpendicular to the direction of the blade. The configuration space is $Q=SE(2)$, with coordinates $(x_1,x_2,\theta)$. The coordinates $(x_1,x_2)$ denote the contact point of the blade with the plane and $\theta$ the orientation of the blade. The Lagrangian is of kinetic type and  if we assume that the center of mass lies in the line through the blade then it is given by
$$
L=\frac{1}{2}\left( (J+ma^2)\dot{\theta}^2+m\left(\dot{x}_1^2+\dot{x}_2^2+2a\dot{\theta}(-\dot{x}_1\sin(\theta)+\dot{x}_2\cos(\theta))\right) \right),
$$
where $m$ denotes the mass of the body, $J$ the moment of inertia relative to the center of mass and $a$ the distance between the center of mass and the contact point of the blade. The matrix of the metric defining the kinetic Lagrangian is then given by
$$
\left(
\begin{array}{ccc}
m & 0 & -ma\sin(\theta) \\
0 & m & ma\cos(\theta) \\
-ma\sin(\theta) & ma\cos(\theta) & J+ma^2
\end{array}
\right) \, .
$$
The nonholonomic constraint is $-\dot{x}_1\sin(\theta)+\dot{x}_2\cos(\theta)=0$, which defines the non-integrable distribution 
$$
{\mathcal D}=\mbox{span}\left\{\frac{\partial}{\partial \theta},  \cos(\theta)\frac{\partial}{\partial x_1}+\sin(\theta)\frac{\partial}{\partial x_2}  \right\} .
$$
For more details on this system, see \cite{neimarkfufaev}.

In \cite{FZ} there is a qualitative study of the DLA method when applied to the Chaplygin sleigh. More precisely, it is shown that the discrete momentum dynamics reproduces the same qualitative behaviour as the continuous momentum dynamics, as long as $\mid \theta'-\theta \mid<2\pi$ and the momentum variable $\rho_2$ satisfies some bound. In the present example we examine the same issue using a discrete gradient method to the equations in adapted coordinates. We will obtain a bound on $h$ but no bound on the momentum variables.

To derive the equations in adapted coordinates, we choose the following orthonormal basis adapted to $D$ and $D^{\bot,g}$ :
\begin{align*}
{\mathcal D}&=\mbox{span}\left\{X_1=\frac{1}{\sqrt{J+ma^2}}\frac{\partial}{\partial \theta}\, , \, X_2=\frac{1}{\sqrt{m}}\left( \cos(\theta)\frac{\partial}{\partial x_1}+\sin(\theta)\frac{\partial}{\partial x_2} \right) \right\} \, ,\\
{\mathcal D}^{\bot,g}&=\mbox{span}\left\{ X_3=\frac{1}{\sqrt{\frac{(J+ma^2)^2}{ma^2}-(J+ma^2)}}\left(\frac{(J+ma^2)}{ma}\sin(\theta)\frac{\partial}{\partial x_1}-\frac{(J+ma^2)}{ma}\cos(\theta)\frac{\partial}{\partial x_2}+\frac{\partial}{\partial \theta}\right) \right\} \, ,
\end{align*}
and denote the induced coordinates on $TQ$ by $(x_1,x_2,\theta,y^1,y^2,y^3)$. In these coordinates the restricted Lagrangian $l:{\mathcal D}\longrightarrow \mathbb{R}$ is given by $l(q^i,y^a)=\frac{1}{2}((y^1)^2+(y^2)^2)$, and the nonholonomic constraint by $y^3=0$.

Since we have chosen an orthonormal basis we have the restricted Hamiltonian $\mathcal{H}(q^i,\rho_1,\rho_2)=\frac{1}{2}\left(\rho_1^2+\rho_2^2\right)$, where $\rho_1=\frac{\partial l}{\partial y^1}=y^1, \rho_2=\frac{\partial l}{\partial y^2}=y^2$. Then
$$
\Pi(\theta,\rho_1)=
\left(
\begin{array}{ccccc}
0 & 0 & 0 & 0 & \frac{\cos(\theta)}{\sqrt{m}} \\
0 & 0 & 0 & 0 & \frac{\sin(\theta)}{\sqrt{m}} \\
0 & 0 & 0 & \frac{1}{\sqrt{J+ma^2}} & 0 \\
0 & 0 & -\frac{1}{\sqrt{J+ma^2}} & 0 & -\frac{a\sqrt{m}}{J+ma^2} \rho_1 \\
-\frac{\cos(\theta)}{\sqrt{m}} & -\frac{\sin(\theta)}{\sqrt{m}} & 0 & \frac{a\sqrt{m}}{J+ma^2} \rho_1 & 0  
\end{array}
\right) \, ,
\quad  \quad 
\nabla \mathcal{H}(\zeta)=\left(
\begin{array}{c}
0 \\
0 \\ 
0 \\ 
\rho_1 \\
\rho_2
\end{array}
\right)\, ,
$$
and the equations of motion \eqref{dynamical} are for the position
\[
\dot{x}_1=\frac{\cos(\theta)}{\sqrt{m}} \rho_2, \quad \dot{x}_2=\frac{\sin(\theta)}{\sqrt{m}} \rho_2, \quad \dot{\theta}=\frac{1}{\sqrt{J+ma^2}}\rho_1\, , \nonumber \\ 
\]
and for the momenta
\begin{equation}
\label{mom-eq}
\dot{\rho_1}=-\frac{a\sqrt{m}}{J+ma^2}\rho_1\rho_2 , \quad \dot{\rho_2}=\frac{a\sqrt{m}}{J+ma^2} \rho_1^2\, .
\end{equation}
The obtained equations are rather simple, since we have a quadratic vector field, a quadratic Hamiltonian and no constraints.

The mean value discrete gradient \eqref{AVF}, the midpoint discrete gradient \eqref{gonzales} and the coordinate increment discrete gradient \eqref{itoAbe} all coincide and give $\bar{\nabla}\mathcal{H}(\zeta,\zeta')=\left(0,0,0,\frac{\rho_1+\rho'_1}{2}, \frac{\rho_2+\rho'_2}{2} \right)^T$. As an approximation to the matrix $\Pi$ we have chosen the midpoint value $\tilde{\Pi}(\zeta,\zeta')=\Pi\left(\frac{\zeta+\zeta'}{2}\right)$. Recalling Remark \ref{rem:IM}, the energy-preserving integrator \eqref{E-int} with any of these discrete gradients then collapses to the implicit midpoint rule, and is consequently given by 
\begin{subequations}
\label{ChapEq}
\begin{align}
x_1'&=x_1+\frac{h}{2\sqrt{m}}\cos \left(\frac{\theta+\theta'}{2}\right)(\rho_2+\rho'_2) \, , \\
x_2&=x_2+\frac{h}{2\sqrt{m}}\sin\left(\frac{\theta+\theta'}{2}\right)(\rho_2+\rho'_2) \, , \\
\theta'&=\theta+\frac{h}{2}\frac{1}{\sqrt{J+ma^2}}(\rho_1+\rho'_1)  \, ,  \\
\rho'_1&=\rho_1-\frac{h}{4}\frac{a\sqrt{m}}{J+ma^2}(\rho_1+\rho'_1)(\rho_2+\rho'_2) \, ,  \label{p1}\\
\rho'_2&=\rho_2+\frac{h}{4}\frac{a\sqrt{m}}{J+ma^2}(\rho_1+\rho'_1)^2 \, . \label{p2}
\end{align}
\end{subequations}

We will write the equations (\ref{p1}) and (\ref{p2}) as
\begin{align*}
F(\rho_1,\rho_2,\rho'_1,\rho'_2,h)&:=\rho'_1-\rho_1+\frac{h}{4}C_{12}^1 (\rho_1+ \rho'_1)(\rho_2+\rho'_2)=0 \, ,\\
G(\rho_1,\rho_2,\rho'_1,\rho'_2,h)&:=\rho'_2-\rho_2+\frac{h}{4}C_{21}^1(\rho_1+\rho'_1)^2=0 \, .
\end{align*}
Notice that $F(\rho_1,\rho_2,\rho_1,\rho_2,0)=G(\rho_1,\rho_2,\rho_1,\rho_2,0)=0$ and compute
$$
\left.
\left(
\begin{array}{cc}
\frac{\partial F}{\partial \rho'_1} & \frac{\partial F}{\partial \rho'_2} \\
\frac{\partial G}{\partial \rho'_1} & \frac{\partial G}{\partial \rho'_2}
\end{array}
\right)
\right|_{(\rho_1,\rho_2,\rho_1,\rho_2,0)}
=
\left.
\left(
\begin{array}{cc}
1+\frac{h}{4}C_{12}^1 (\rho_2+\rho'_2) & \frac{h}{4}C_{12}^1 (\rho_1+\rho'_1) \\
\frac{h}{2}C_{21}^1 (\rho_1+\rho'_1) & 1
\end{array}
\right)
\right|_{(\rho_1,\rho_2,\rho_1,\rho_2,0)}
=
\left(
\begin{array}{cc}
1 & 0 \\
0 & 1
\end{array}
\right) \, .
$$
By the implicit function theorem we can write $\rho'_1=f(\rho_1,\rho_2,h)$ and $\rho'_2=g(\rho_1,\rho_2,h)$ in a neighbourhood of $(\rho_1,\rho_2,0)$, with $(\rho'_1,\rho'_2)$ also in a neighbourhood of $(\rho_1,\rho_2)$.

The continuous system has certain qualitative characteristics. Specifically, for the continuous system we have from \eqref{mom-eq}, in the case $a\not=0$, a one-dimensional manifold of equilibria $\left\{ \rho_1=0 \right\}$. These equilibria are stable and asymptotically stable with respect to $\rho_1$ if $\rho_2>0$ and unstable if $\rho_2<0$. We will now study how the qualitative behaviour of (\ref{p1})-(\ref{p2}) compares, as in \cite{FZ}. 

\textit{Equilibria}:
If $h\not=0$ then $F(\rho_1,\rho_2,\rho_1,\rho_2,h)=0$ and  $G(\rho_1,\rho_2,\rho_1,\rho_2,h)=0$ imply $\rho_1=0$. Then the set $\left\{ \rho_1=0 \right\}$ is a one-dimensional manifold of equilibria.

\begin{comment}
Notice that if $\rho_1=0$, then $\rho'_1 (1+\frac{h}{4}C_{12}^1(\rho_2+\rho'_2))=0$ and in a sufficiently small neighbourhood of $(\rho_1,\rho_2,\rho_1,\rho_2,0)$ the only solution is $\rho'_1=0$, and therefore $\rho_2=\rho'_2$. 

Notice also that for $\rho_2>0$ the solution is defined for any $h$ while for $\rho_2<0$ we need to require a small enough time step, namely $h<\frac{-2}{C_{12}^1 \rho_2}$.
\end{comment}

\textit{Stability}:
Now we study the linearization of $(f,g)$ at the equilibrium points $eq=(0,\rho_2,0,\rho_2,h)$. Assuming that $\rho_2\not=0$ and $h<\left|\frac{2}{C_{12}^1 \rho_2}\right|$ we compute
$$
\left.
\left(
\begin{array}{cc}
\frac{\partial f}{\partial \rho_1} & \frac{\partial f}{\partial \rho_2} \\
\frac{\partial g}{\partial \rho_1} & \frac{\partial g}{\partial \rho_2}
\end{array}
\right)\right|_{eq}
=
\left(
\begin{array}{cc}
\frac{2-hC_{12}^1\rho_2}{2+hC_{12}^1\rho_2} & 0 \\
0 & 1
\end{array}
\right) \, ,
$$
with eigenvalues
$$
\lambda_1=\frac{2-hC_{12}^1\rho_2}{2+hC_{12}^1\rho_2}=\frac{2(J+ma^2)-ha\sqrt{m}\rho_2}{2(J+ma^2)+ha\sqrt{m}\rho_2} \, , \quad \lambda_2=1 \, .
$$
Since $0<h<\left|\frac{2(J+ma^2)}{a\sqrt{m}\rho_2}\right|=\left|\frac{2}{C_{12}^1\rho_2}\right|$, we have $\lambda_1>0$, regardless of $\rho_2\not=0$. Further if $\rho_2>0$ then $\lambda_1<1$ and hence the equilibrium is stable and asymptotically stable with respect to $\rho_1$. On the other hand if $\rho_2<0$ then  $\lambda_1>1$ and hence the equilibrium is unstable. Therefore the proposed discrete method reproduces the same qualitative behaviour as the continuous system. This is not guaranteed when applying the midpoint rule to the Chaplygin sleigh system in the original coordinates.

\begin{prop} \label{cs-stability}
The energy-preserving method (\ref{p1})-(\ref{p2}) has a one-dimensional manifold of equilibria $\left\{ \rho_1=0 \right\}$. Assuming that $h<\left|\frac{2}{C_{12}^1 \rho_2}\right|$, the equilibria $(0,\rho_2)$ are stable and asymptotically stable with respect to $\rho_1$ if $\rho_2>0$ and are unstable if $\rho_2<0$. 
\end{prop}

\begin{minipage}[t]{0.45\textwidth}
\centering
   \includegraphics[width=\textwidth]{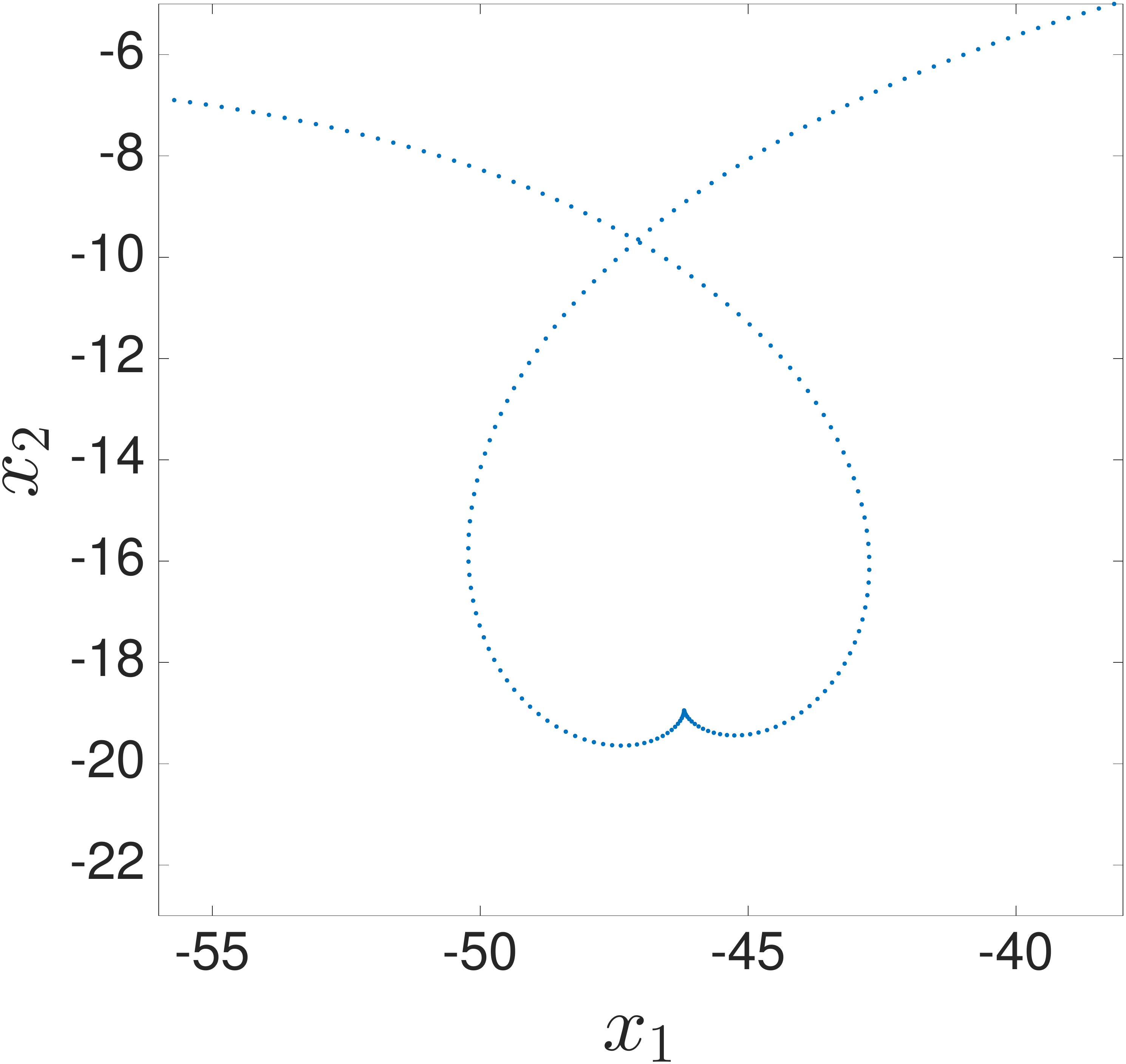}
 
\end{minipage}
\hspace{0.05\linewidth}
\begin{minipage}[t]{0.45\textwidth}
\centering
   \includegraphics[width=\textwidth]{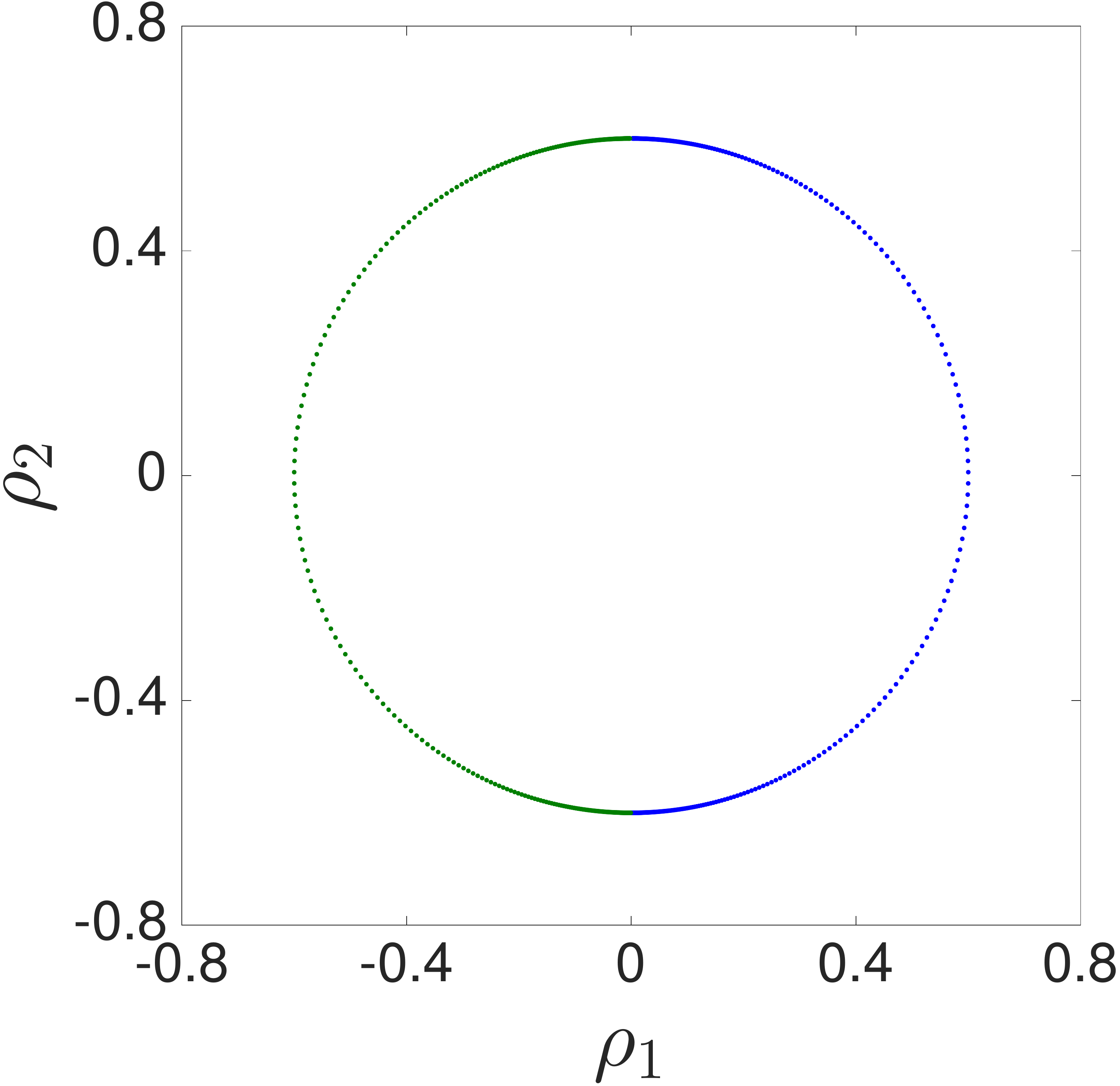}
\end{minipage}
\begingroup
 \captionof{figure}{{\footnotesize Integration results for the sleigh, using the method \eqref{ChapEq}, with the parameters set to $J=8$, $a=m=1$, step-size $h=0.5$ and initial values $x_1=-5,\, x_2=0,\, \theta = 0.1,\,\rho_1\in\left\{-0.001,0.001\right\},\, \rho_2=-0.6$. {\bf Left}: A partial $x_1x_2$ trajectory. {\bf Right}: Two $\rho_1 \rho_2$ trajectories.} \label{fig:Chap}}
\endgroup
\medskip

In Figure \ref{fig:Chap} we see an example of how the method exhibits correct behaviour by converging towards a stable equilibrium point when starting very close to an unstable one.
\vspace{0.3cm}
\begin{remark}
Similarly it is possible to show that any convergent Runge-Kutta method will give the correct behaviour for small enough $h$, when applied to the equations (\ref{mom-eq}). For example, applying the explicit Euler method to these equations, we obtain the same conclusion as in Proposition \ref{cs-stability} if we assume $h<\left|\frac{1}{C_{12}^1 \rho_2}\right|$. This confirms the fact that the illustrated good qualitative behaviour with respect to stability of equilibria is more an effect of the choice of coordinates rather than of the choice of the method applied in those coordinates.
\end{remark}

\subsection{Euler-Poincar\'{e}-Suslov problem on $\mathfrak{so}(3)$}\label{section43}
In this example we show that the approach that we have presented is also valid for nonholonomic systems defined on a Lie algebra (and more generally on a Lie algebroid \cite{CLMM09}). 

Let $\left\{ e_1,e_2,e_3 \right\}$ be a basis of the Lie algebra $\mathfrak{so}(3)\cong {\mathbb R}^3$ and denote the corresponding  coordinates  by $(\omega_1,\omega_2,\omega_3)$. Consider a kinetic Lagrangian on $\mathfrak{so}(3)$  defined by the matrix
$$
(\textbf{\textit{g}}_{ij}) = \left(
\begin{array}{ccc}
I_{11} & I_{12} & I_{13} \\
I_{12} & I_{22} & I_{23} \\
I_{13} & I_{23} & I_{33}
\end{array}
\right) \, ,
$$
and introduce the nonholonomic constraints given by $\sum a_i\omega_i=0$, where $a\in\mathfrak{so}(3)$ is a fixed element. We can choose the frame $\left\{ e_1,e_2,e_3 \right\}$ in such a way that $I_{12}=0$ and $a=e_3$. Then the Lagrangian is given by
$$
L=\frac{1}{2}(I_{11}\omega_1^2+I_{22}\omega_2^2+I_{33}\omega_3^2+2I_{13}\omega_1\omega_3+2I_{23}\omega_2\omega_3)\, ,
$$
and the constraint reduces to $\omega_3=0$. This defines the distribution 
$$
\mathcal{D}=\mbox{span}\left\{ X_1:=(1,0,0),\, X_2:=(0,1,0) \right\}.
$$
Since the bracket on $\mathfrak{so}(3)$ is given by the cross product, it is immediate that $\mathcal{D}$ is not involutive, and hence the constraint $\omega_3=0$ is nonholonomic.
On the other hand,
$$
\mathcal{D}^{\bot,g}=\mbox{span}\left\{ X_3:=(I_{22}I_{13}, I_{11}I_{23}, -I_{11}I_{22}) \right\} \, .
$$
The Lie bracket of $X_1$ and $X_2$ is expressed in terms of the adapted basis $X_1$, $X_2$, $X_3$ as
$$
\left[ X_1, X_2 \right]=(1,0,0)\times (0,1,0)=(0,0,1)=\frac{I_{13}}{I_{11}}X_1+\frac{I_{23}}{I_{22}}X_2-\frac{1}{I_{11}I_{22}}X_3 \, .
$$
Thus the nonvanishing structure constants of the projected bracket are
$$
C_{12}^1=\frac{I_{13}}{I_{11}} \quad \mbox{and} \quad C_{12}^2=\frac{I_{23}}{I_{22}}.
$$
If we denote by $(y^1,y^2,y^3)$ the coordinates corresponding to the adapted basis $\left\{X_1, X_2, X_3\right\}$, the change of coordinates is given by 
$$
\omega_1=y^1+I_{22}I_{13}y^3 \, , \quad  \omega_2=y^2+I_{11}I_{23}y^3 \, , \quad  \omega_3=-I_{11}I_{22}y^3 \, .
$$ 
Then the restricted Lagrangian becomes $l=\frac{1}{2}(I_{11}(y^1)^2+I_{22}(y^2)^2)$ and the nonholonomic constraint is $y^3=0$.

In this example, since there are no $(q^i)$ variables, the equations of motion (\ref{H1})-(\ref{H2}) reduce to $\dot{\rho}_a=-C_{ab}^c \rho_c\frac{\partial \mathcal{H}}{\partial \rho_b}$, where $\mathcal{H}=\frac{1}{2}\left(\frac{1}{I_{11}}\rho_1^2 + \frac{1}{I_{22}}\rho_2^2 \right)$ and $\rho_i=\frac{\partial l}{\partial y^i}=I_{ii}y^i$, that is
$$
\dot{\rho}_{1}=-\frac{1}{I_{22}}(C_{12}^1\rho_1+C_{12}^2\rho_2)\rho_2 \quad \mbox{and} \quad  \dot{\rho}_{2}=-\frac{1}{I_{11}}(C_{21}^1\rho_1+C_{21}^2\rho_2)\rho_1 \, .
$$

In matrix form, using $C_{ab}^c = -C_{ba}^c$, we get
$$
\left(
\begin{array}{c}
\dot{\rho}_{1} \\
\dot{\rho}_{2}
\end{array}
\right)
=
-\left(
\begin{array}{cc}
0 & C_{12}^1\rho_1+C_{12}^2\rho_2 \\
-(C_{12}^1\rho_1+C_{12}^2\rho_2) & 0
\end{array}
\right)
\left(
\begin{array}{c}
\frac{\rho_1}{I_{11}} \\
\frac{\rho_2}{I_{22}}
\end{array}
\right) \, .
$$

We apply the same discrete gradient method as in the previous example to get the integrator
$$
\left(
\begin{array}{c}
\frac{\rho'_1-\rho_1}{h} \\
\frac{\rho'_2-\rho_2}{h}
\end{array}
\right)
=
-\left(
\begin{array}{cc}
0 & C_{12}^1\frac{\rho'_1+\rho_1}{2}+C_{12}^2\frac{\rho'_2+\rho_2}{2} \\
-\lp C_{12}^1\frac{\rho'_1+\rho_1}{2}+C_{12}^2\frac{\rho'_2+\rho_2}{2}\rp & 0
\end{array}
\right)
\left(
\begin{array}{c}
\frac{\rho'_1+\rho_1}{2I_{11}} \\
\frac{\rho'_2+\rho_2}{2I_{22}}
\end{array}
\right) \, .
$$
As mentioned previously the integrator is here equivalent to the implicit midpoint method.

\subsection{Continuous gearbox driven by an asymmetric pendulum}\label{section44}

In this final example we compare the performance of integrators applied directly to the formulation of the system in canonical coordinates.  We consider a continuous gearbox driven by an asymmetric pendulum. This is a special case of the continuous gearbox system discussed in \cite{MV14}.  Here $Q = \mathbb{R}^3$ with Hamiltonian
\begin{subequations}
\label{asyPend}
\begin{align}
H(q^i,p_i) &= \frac{1}{2}\left(p_1^2+p_2^2+p_3^2\right)-V(q^i) \, , \\  
V(q^i) &= \frac{1}{2}\left((q^1)^2+(q^2)^2\right)+\cos(q^3)-\frac{1}{5}\sin(2q^3)\, .
\end{align}
\end{subequations}
The single nonholonomic linear velocity constraint is 
\[
\dot{q}^1+\sin(q^3)\dot{q}^2 = p_1+\sin(q^3)p_2 = 0 \, ,
\]
since clearly $p_i = \dot{q}^i$, $1 \leq i \leq 3$ because $(\textbf{\textit{g}}_{ij}) = I_3$. 

Again comparing with the semi implicit reversible DLA variational integrator (SI-DLA) proposed in \cite{perlmutter06}, we consider from Section \ref{T*Q} the initial method \eqref{simpDirect} using the Gonzalez midpoint discrete gradient (Gonzalez-F) and the canonical coordinate implementation of the Gonzalez midpoint discrete gradient method for the reduced system (Gonzalez-R). For the methods we compare in Figure \ref{fig:asypend} the relative energy error, i.e. $\left|(H-H_0)/H_0\right|$, and the constraint error, i.e. $|p_1+\sin(q^3)p_2|$ for a long time simulation $t\in [0,50000]$ with random initial values chosen to ensure non-periodic behaviour. $H_0$ is the initial energy.

For SI-DLA we observe an exponential growth in the energy error, while both Gonzalez methods preserve the energy to machine precision. SI-DLA and GONZALEZ-R both preserve the nonholonomic constraint to machine precision. However the energy blow-up of SI-DLA gives a corresponding increase in the round off error and thus the constraint error over time. GONZALEZ-F does not respect the nonholonomic constraint. The results are thus as expected.

\medskip

\begin{minipage}[t]{0.45\textwidth}
\centering
\includegraphics[width=\textwidth]{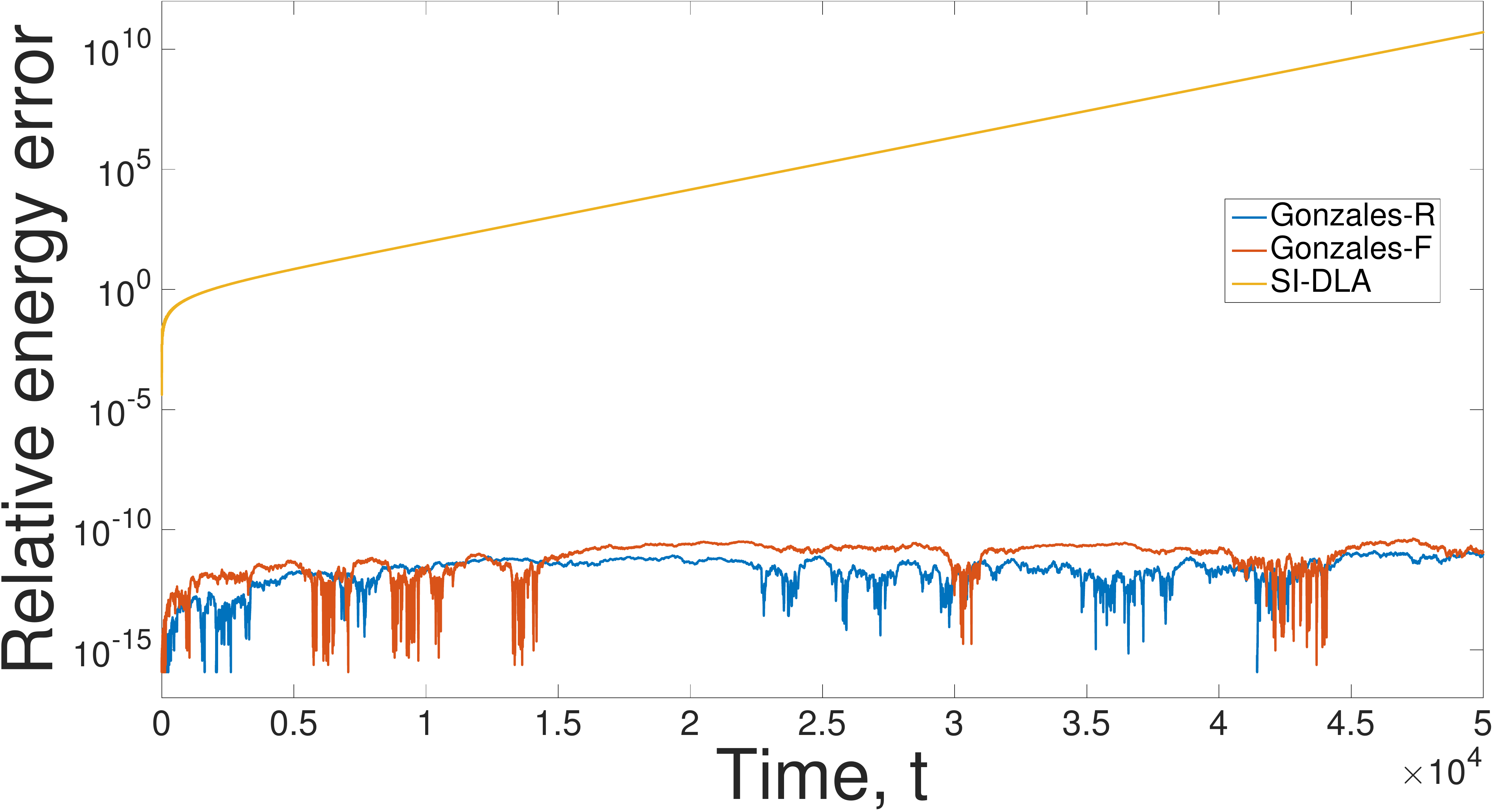}
\end{minipage}
\hspace{0.05\linewidth}
\begin{minipage}[t]{0.45\textwidth}
\centering
\includegraphics[width=\textwidth]{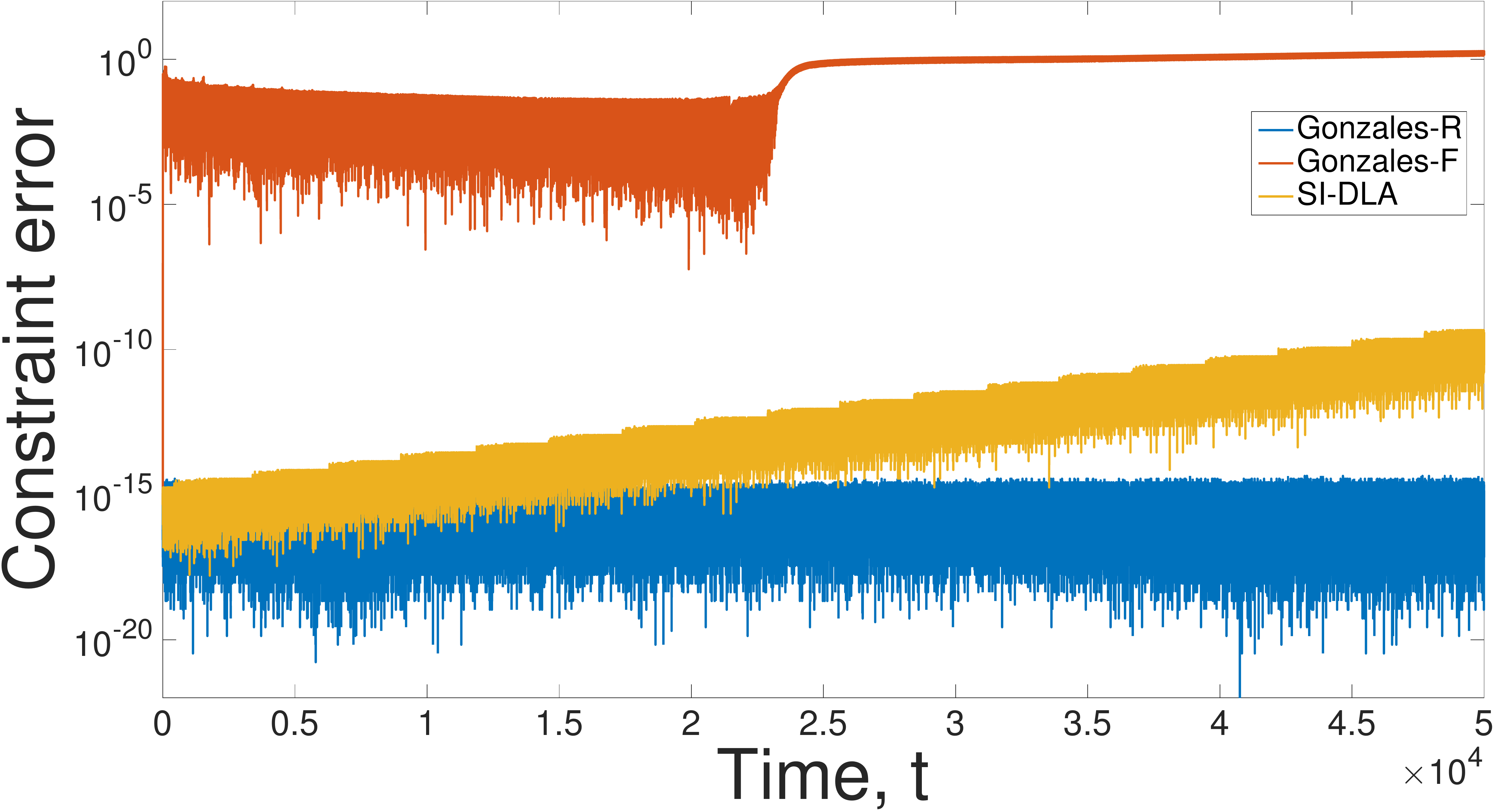}
\end{minipage}
\begingroup
\captionof{figure}{ {\footnotesize Integration results for methods GONZALEZ-R, GONZALEZ-F and SI-DLA, applied to the system \eqref{asyPend} with random initial values, $h=0.1$, and $t \in [0,50000]$. {\bf Left}: Relative energy error, i.e. $|H(t)-H(0)/H(0)|$. {\bf Right} The error in the nonholonomic constraint, i.e. $|p_1+\sin(q^3)p_2|$. \label{fig:asypend}}}
\endgroup

\section{Future work}
 In this paper we have not  addressed the case  when  $Q$ is a differentiable manifold. In a future paper, we will  propose to adapt the discrete gradient approach taking the geometry of the configuration space into account, see for instance the methods in reference \cite{CeOw}. In order to adapt the ideas in \cite{CeOw} to a general differentiable manifold $Q$, we will need to introduce 
a finite difference map  or retraction map $\Phi_h: U\subset {\mathcal D}^*\times {\mathcal D}^*\rightarrow T{\mathcal D}^*$  (see  \cite{perlmutter06}) from a finite difference map initially defined on $Q$.  In this case we will define a discrete gradient as a map
$\bar{\nabla}{\mathcal H}:{\mathcal D}^*\times {\mathcal D}^* \longrightarrow T^{*}{\mathcal D}^*$ verifying similar properties to Definition \ref{def31} (see \cite{CeOw}).
$$
\xymatrix{
{\mathcal D}^*\times {\mathcal D}^* \ar[rr]^{\bar{\nabla}{\mathcal H}} \ar[d]^{\Phi_h} && T^{*}{\mathcal D}^* \ar[d]^{\pi_{{\mathcal D}^*}}\\
	 T{\mathcal D}^*  \ar[rr]^{\tau_{{\mathcal D}^*}} && {\mathcal D}^* 
		}
		$$
In this case, an energy preserving integrator  for Equation (\ref{dynamical}) would be 
\[
 \Phi_h(\zeta, \zeta')=\Pi(\bar{\zeta})\bar{\nabla}{\mathcal H}(\zeta, \zeta')
\]
with $\bar{\zeta}=\tau_{\mathcal D}^*(\Phi_h(\zeta, \zeta'))$. 
We will explore this possibility in a future paper since in many examples of nonholonomic systems the configuration space is a nonlinear space such as, for instance, a Lie group $G$. 

Moreover, it would be interesting to compare the discrete gradient method approach introduced in this paper with other methods designed for nonholonomic systems. For instance, the Chaplygin case is given by  a Lagrangian system with forces on the tangent space of a reduced space and then it is possible to use directly discrete variational integrators based on forced Lagrangian systems (see \cite{cortmart,cortes02}). Other interesting possibilities to compare our methods with are variational integrators from Hamiltonizable nonholonomic
              systems \cite{FeBoOl} or the geometric nonholonomic integrator \cite{FeIgMa}.

\section*{Acknowledgments}

This work has been partially supported by MINECO (Spain) MTM2013-42870-P, MTM2015-69124-REDT, the ICMAT Severo Ochoa project SEV-2015-0554 and the Nils-Abel project
010-ABEL-CM-2014ANILS. It has also received funding from the European Unions
Horizon 2020 research and innovation programme under the Marie Sklodowska-Curie
grant agreement No. 691070. MFP has been financially supported by a FPU scholarship from MECD.

\bibliographystyle{plain}
\bibliography{References}

\end{document}